\documentclass[11pt]{amsart}
\usepackage{amsmath,amssymb,mathrsfs,amscd}
\usepackage{graphicx}

\title{On a construction of the twistor spaces of Joyce metrics, I}
\author{Nobuhiro Honda
}  

\newcommand{\ol}{\overline}
\newcommand{\ra}{\rightarrow}
\newcommand{\lra}{\longrightarrow}

\newcommand{\Ra}{\Rightarrow}
\newcommand{\set}{\,|\,}
\newcommand{\proofend}{\hfill$\square$}

\newcommand{\vsp}{\vspace{3mm}}

\newtheorem{prop}{Proposition}[section]
\newtheorem{lemma}[prop]{Lemma}

\newtheorem{thm}[prop]{Theorem}
\newtheorem{rmk}[prop]{Remark}

\begin{document}

\maketitle

\begin{abstract} 
We  explicitly construct  the twistor spaces of some self-dual metrics with torus action given by D.\,Joyce. 
Starting from a fiber space over a projective line whose fibers are compact singular toric surfaces, we apply a number of birational transformations to obtain the desired twistor spaces. 
Especially an important role is played by Atiyah's flop.
\end{abstract}



\section{Introduction} 
In 1995, D.\,Joyce \cite{J95} constructed a series of self-dual metrics with torus action on the connected sum of complex projective planes.
A.\,Fujiki  \cite{F00} intensively studied algebro-geometric properties of the twistor spaces of Joyce metrics and obtained a characterization of Joyce metrics saying that  any self-dual metric on these 4-manifolds admitting an effective torus action must coincide with  Joyce metrics.
However, it is still unknown the way how to construct the twistor spaces.
In this paper we shall explicitly construct the twistor spaces of some Joyce metrics on $4\mathbf P^2$.


Unlike past algebraic studies  of twistor spaces, the main tool in this paper is the anticanonical system of the twistor spaces.
In Section 2, 
after  reviewing some  necessary results in \S 2.1, 
we first show in \S 2.2 that the meromorphic map associated to the anticanonical system of the twistor space is  bimeromorphic onto its image.
Then next we  determine its defining equations of the image in a projective space.
The conclusion is that the anticanonical model of the twistor space is a complete intersection of degree 8 in $\mathbf P^6$ of a particular form (Theorem \ref{thm-01}).
Subsequently in \S 2.3 we obtain more necessary conditions on  the defining equations (Theorem \ref{thm-02}),  concerning how the defining quadratics intersect. 

Section 3 (and the figures there) is a heart of this paper. Based on the computations in Section 2,  we will give an explicit way how one can obtain the twistor space starting from the projective model. 
Our procedure (Step 1-- Step 6) consists of a number of blowing-ups and blowing-downs, including careful choice of small resolutions of ordinary double points, as well as many flops which are always possible if there is a smooth rational curve in a 3-fold having an appropriate normal bundle.
All the operations are concretely displayed in figures.

We remark that a reader who is interested in the construction itself can consult Section 3 and Figures \ref{fig-bim0}--\ref{fig-bim6}  directly. 
The results of Section 2 are needed in order to prove that the constructed 3-folds are actually the twistor spaces of Joyce metrics.

It should also be worth mentioning that the defining equations of the anticanononical models of the twistor spaces  are obtained in effect  in the following manner:
(i) first  the defining equations of the bi-anticanonical model of a toric surface (in $\mathbf P^4$) are determined,
(ii) then the coefficients of the defining equations are replaced by polynomials, and they are thought as equations defined on a ($\mathbf P^4$-)bundle over a projective line.


\vsp\noindent
{\bf{Notations and Conventions.}}
Throughout this paper we try to use the same notations as those in \cite{F00}.
In particular, $K=U(1)\times U(1)$,  $G=\mathbf C^*\times\mathbf C^*$ and
for a positive integer $n$, $n\mathbf{P}^2$ denotes the connected sum of $n$ copies of complex projective planes.
For a complex vector space $V$, $\mathbf P^{\vee}V$ denotes the dual complex projective space.
The dimension of a linear system always means the projective dimension, so that if $|L|$ is the complete linear system associated to  a holomorphic line bundle $L$ over a compact complex manifold $X$, the associated map is a meromorphic map from $X$ to $\mathbf P^{\vee}H^0(X,L)=\mathbf P^{n-1}$, where $n=\dim H^0(X,L)$.
Further, Bs\,$|L|$ denotes the base locus of the linear system.
For a linear subspace $P$ of a projective space, $\dim P=0$ implies that $P$ is a point.
A smooth rational curve in a 3-fold is called a $(-1,-1)$-curve if the 3-fold is smooth in a neighborhood of the curve and if its normal bundle is isomorphic to $\mathscr O(-1)\oplus \mathscr O(-1)$.
If $A$ is a finite set, $\#A$ denotes the number of elements of $A$.

\section{The anticanonical models of the twistor spaces}
\subsection{Review of known results on the twistor spaces of Joyce metrics}
First of all we briefly recall a description of $K$-actions on $n\mathbf{P}^2$ following \cite{F00, J95} in order to explain which Joyce metrics we are concerned.
The set of equivalence classes of all smooth effective $K$-actions on $n \mathbf{P}^2$ can be classified by using combinatorial data, in an analogous way to a description of toric surfaces.
Concretely, for a given smooth effective $K$-action on $n\mathbf{P}^2$, the singular set of the action, namely the set of points having a non-trivial isotropy group, forms a `cycle' of  $(n+2)$ spheres. The intersection points of these spheres are fixed points of the $K$-action and if we remove these points from the cycle of  singular set, the isotropy subgroup of each connected component becomes a $U(1)$-subgroup
$$
K(l, m)=\{(s,t)\in K\set s^lt^m=1 \},
$$
where $l$ and $m$ are relatively coprime integers determined only up to similtaneous inversion of signs.
If $K(l,m)$ and $K(l',m')$ are the isotropy subgroups of adjacent components, we have
$
|lm'-ml'|=1.
$
Conversely, if we are given  a sequence of $(n+2)$ subgroups of the form $K(l,m)$ 
satisfying these relations,  one can construct a smooth effective $K$-action on $n\mathbf{P}^2$ having the sequence as the isotropy subgroups of the singular set.

For each  smooth effective  $K$-action on $n\mathbf P^2$, Joyce \cite{J95} constructed an $(n-1)$-dimensional family of self-dual metrics on $n\mathbf{P}^2$ having the given $K$-action as an isometry group.
One can also consult a paper of Pedersen-Poon \cite{PP95}, where the existence of such self-dual metrics is proved by using a method of Donaldson and Friedman.
Conversely,  Fujiki \cite{F00} proved that any self-dual metric on $n\mathbf P^2$ admitting a smooth effective $K$-action is isomorphic to a Joyce metric.
When $1\le n\le 3$, Joyce metrics on $n\mathbf P^2$ are not new in the sense that they are nothing but a special form of LeBrun metrics \cite{LB91}.
This is partly because   smooth effective $K$-actions on $n\mathbf P^2$ is unique up to equivalence, for $1\le n\le 3$, according to a formula of Fujiki \cite[Prop.\,\,3.3]{F00}.
The formula also shows that if $n=4$, there arise 3 kinds of $K$-actions that are mutually non-equivalent.
Therefore on $4\mathbf P^2$ there are precisely 2 families of Joyce metrics that are different from LeBrun metrics.
These  2 kinds of $K$-actions respectively have the following  circular sequences of $U(1)$-subgroups as their isotropy:

$$\{ K(1,0),K(1,1), K(0,1), K(-1,2), K(-2,3),K(-1,1)\}\,\cdots\text{Type I},$$
$$\{K(1,0),K(1,1),K(0,1),K(-1,2),K(-1,1), K(-2,1)\}\,\cdots\text{Type II}.$$

In this paper we are interested in the twistor spaces of Joyce metrics on $4\mathbf P^2$ whose $K$-action is Type I.
Before starting its investigation, we recall results of Fujiki \cite{F00} about  beautiful structures of the twistor spaces of  Joyce metrics on $n\mathbf P^2$ in general.
For this,  recall that if a torus $K$ acts on a compact 4-manifold preserving a self-dual structure, then  the $K$-action naturally induces a holomorphic $G$-action on $Z$. 
We also recall that if $S$ is a compact smooth toric surface, and if $T\,(\simeq G)$ is the unique 2-dimensional orbit in $S$, then $C=S\backslash T$ is the unique $G$-invariant anticanonical curve on $S$, which is necessarily a cycle of $k$ rational curves where $k=12-K^2_S$.

\begin{prop}\label{prop-01}\cite{F00}
Let $Z$ be the twistor space of a Joyce metric on $n\mathbf P^2$ which is different from LeBrun metrics. Then  the following holds.
(i) $\dim |(-1/2)K_Z|=1$ and its members are  $G$-invariant.
Moreover, general members of the pencil are irreducible and smooth, and they become a toric surface with $K^2=8-2n$.
(ii) Let $S$ be a general member of the pencil 
and  $C$ the  unique $G$-invariant anticanonical curve on $S$. 
Then the base locus of the pencil is $C$.
(iii) If we write $C=\sum_{i=1}^{2n+4}C_i$ in such a way that $C_i$ and $C_{i+1}$ intersect, the real structure on $Z$ exchanges $C_i$ and $C_{i+n+2}$, where the subscriptions are counted modulo $2n+4$.
(iv) The pencil $|(-1/2)K_Z|$ has precisely $(n+2)$ singular members and all of them consist of two smooth irreducible toric surfaces that are conjugate of each other.
(v) If we write $S_i=S_i^++S_i^-$ ($1\le i\le n+2$) for the reducible members, $S^+_i$ and $S^-_i$ divide $C$ into 'halves' in the sense that both $S_i^+\cap C$ and $S_i^-\cap C$ are connected (cf.\! (iii)). 
Moreover,  $L_i:=S^+_i\cap S_i^-$ is a $G$-invariant twistor line. 
\end{prop}

The restriction of the twistor fibration $Z\ra n\mathbf P^2$ on a smooth member $S$ of $|(-1/2)K_Z|$ is $K$-equivariant and we can determine the structure of $S$ as a toric surface from the $K$-action on $n\mathbf P^2$, by noting that the restriction of the twistor fibration on $C$ gives an unramified double covering onto the singular set of the $K$-action on $n\mathbf P^2$.
Also, it is possible to determine the structure of $S_i^+$ (and $S_i^-$) as toric surfaces for every $1\le i\le n+2$ in a similar way.

\subsection{The anticanonical systems of the twistor spaces}
For the rest of this paper $Z$ always denotes the twistor space of a Joyce metric on $4\mathbf P^2$ whose $K$-action is of Type I.
The purpose of this subsection is  to prove the following result about the anticanonical system  of the twistor space $Z$ and its associated map:

\begin{thm}\label{thm-01}
We have the following.
(i)  $\dim |-K_Z|=6$ and\, {\rm{Bs}}\,$|-K_Z|\subset C$, where $C$ is a cycle of 12 rational curves as in Proposition \ref{prop-01}
(ii) The meromorphic map $\Phi:Z\ra\mathbf{P}^6$ associated to $|-K_Z|$ is bimeromorphic onto its image. (iii) The image of\, $\Phi$ is  a complete intersection in\, $\mathbf P^6$ defined by the following three equations
\begin{equation}
\label{001}
x_1x_2=Q_1(x_5,x_6,x_7),\,\,
x_3x_4=Q_2(x_5,x_6,x_7),\,\,
x_5^2=x_6x_7,
\end{equation}
where $(x_1,\cdots,x_7)$ is a homogeneous coordinate on\, $\mathbf{P}^6$, and both $Q_1$ and $Q_2$ are homogeneous quadratic polynomials of $(x_5,x_6,x_7)$ with real coefficients.
(iv) In the above coordinate the real structure on $\mathbf P^6$ induced from that on $Z$ is given by
\begin{equation}\label{0011}
(x_1,x_2,x_3,x_4,x_5,x_6,x_7)\longmapsto 
(\ol{x}_2,\ol{x}_1,\ol{x}_4,\ol{x}_3,\ol{x}_5,\ol{x}_6,\ol{x}_7).
\end{equation}
\end{thm}

For a proof of Theorem \ref{thm-01} we study bi-anticanonical system of a general member $S$ of the pencil $|(-1/2)K_Z|$.
By Proposition \ref{prop-01}, $S$  is a smooth toric surface with $K^2=0$.
For the cycle $C$ of $G$-invariant anticanonical curve on $S$, we respect the real structure and write
$$
C=\sum_{i=1}^6C_i+\sum_{i=1}^6\ol{C}_i,
$$
where $\ol{C}_i$ denotes the conjugation of $C_i$ (namely the image of $C_i$ by the real structure).
We can suppose
 $C_iC_{i+1}= \ol{C}_i\ol{C}_{i+1}=1$ for $1\le i\le 5$ and $C_6\ol{C}_1=\ol{C}_6 C_1=1$ hold (cf. Figure \ref{fig-6tls}).
Moreover, the self-intersection numbers of $C_i$ inside $S$ can be easily computed and after a possible renumbering we have
$$
C_1^2=-3,\,C_2^2=-2,\,C_3^2=-1,\,C_4^2=-3,\,C_5^2=-2,\,C_6^2=-1.
$$
Then we have the following two results on the bi-anticanonical system of $S$:

\begin{lemma}\label{lemma-01}
Let $S$ be an irreducible member of the pencil $|(-1/2)K_Z|$  and $C=\sum C_i+\sum \ol{C}_i$ the $G$-invariant anticanonical curve as above.
Then we have the following.
(i) {\rm{Bs}}\,$|-2K_S|=C_1+C_2+C_4+C_5+\ol{C}_1+\ol{C}_2+\ol{C}_4+\ol{C}_5$ and the movable part of the system  is base point free and 4-dimensional.
(ii) The associated morphism $\phi:S\ra\mathbf P^4$ is birational onto its image.
(iii) The image of $\phi$ is a complete intersection in $\mathbf{P}^4$ defined by
\begin{equation}\label{002}
x_1x_2=ax_0^2,\,\,x_3x_4=bx_0^2
\end{equation}
for some $a,b\in\mathbf C^*$,
where  $(x_0,\cdots,x_4)$ is a homogeneous coordinate on $\mathbf{P}^4$.
(iv) If $S$ is a real member, the real structure on $\mathbf P^4=\mathbf P^{\vee}H^0(-2K_S)$ induced from $S$ is given by, in the above coordinate, 
\begin{equation}\label{0031}
(x_0,x_1,x_2,x_3,x_4)\longmapsto(\ol{x}_0,\ol{x}_2,\ol{x}_1,\ol{x}_4,\ol{x}_3),
\end{equation}
so that $a$ and $b$ must be real numbers.
\end{lemma}

\begin{lemma}\label{lemma-011} 
(i) If $ab\neq 0$, the surface \eqref{002} in $\mathbf P^4$ is a toric surface whose singular locus consists of  4 ordinary double points.
(ii) If $(a,b)\in\mathbf R^2$, the surface \eqref{002} has a real point iff $a\ge 0$ and $b\ge0$.
(iii) If  $\phi:S\ra\mathbf P^4$ is as in Lemma \ref{lemma-01}, $\phi:S\ra \phi(S)$ is the composition of the blowing-down of the four $(-1)$  curves $C_3, C_6,\ol{C}_3,\ol{C}_6$ and  the contraction of the images of $C_1,C_4,\ol{C}_1,\ol{C}_4$ (which become $(-2)$-curves) to the  ordinary double points of $\phi(S)$.
\end{lemma}

The equations \eqref{002}  will become important in our description of the anticanonical models of the twistor spaces. 
Actually they will be obtained in effect by viewing \eqref{002} as  equations in some $\mathbf P^4$-bundle over $\mathbf P^1$ and allowing $a$ and $b$ to be polynomials. 

\vsp
\noindent
Proof of Lemma \ref{lemma-01}.
By calculating intersection numbers, it can be readily seen that the the fixed component of the system $|-2K_S|$ is  $(C_1+C_2+C_4+C_5)+(\ol{C}_1+\ol{C}_2+\ol{C}_4+\ol{C}_5)$.
Removing this from $-2K_S=2C$, we obtain that the movable part of $|-2K_S|$ is the system $|F|$, where $F=-K_S+C_3+C_6+\ol{C}_3+\ol{C}_6$.
It can be shown by a standard argument that $|F|$ is free, 4-dimensional, and $F^2=4$.
Thus we obtain (i).
It can also be easily seen  that a general member  $A$ of $|F|$  is a smooth elliptic curve.
Since $S$ is rational, the restriction map $H^0(S,F)\ra H^0(A, F|_A)$ is surjective.
Because $\deg F|_A=F^2=4$, $F|_A$ is very ample. 
It follows that $\phi|_A$ is biholomorphic onto its image.
Moreover, by the definition of the map $\phi$, we have $\phi^{-1}(\phi(A))=A$.
Therefore $\phi$ is generically one-to-one and hence birational onto its image.
Hence we obtain (ii).
Next we determine the image $\phi(S)$ in $\mathbf{P}^4$ by using $G$-action.
First since $\phi$ is a birational morphism and $F^2=4$, $\phi(S)$ is a non-degenerate surface of degree 4.
On the other hand, by the theory of toric variety, the natural $G$-action on the cohomology group $H^0(-mK_S)$  can be computed by solving a system of linear inequalities over integers \cite{Oda}.
The result of the computation is that, if $m=2$, the $G$-action takes a simple form
\begin{equation}
\label{now}
(y_0,y_1,y_2,y_3,y_4)\longmapsto (y_0, s^{-1}y_1,sy_2,t^{-1}y_3,ty_4)\hspace{2mm}\text{for}\hspace{2mm}(s,t)\in G
\end{equation}
for a coordinate $(y_0,\cdots,y_4)$ on $H^0(-2K_S)$.
Letting $(x_0,\cdots,x_4)$ be the dual coordinate of $(y_0,\cdots,y_4)$ on $H^0(-2K_S)^{\vee}$, the natural $G$-action on $H^0(-2K_S)^{\vee}$ is given by 
\begin{equation}
\label{003}
(x_0,x_1,x_2,x_3,x_4)\longmapsto (x_0, sx_1,s^{-1}x_2,tx_3,t^{-1}x
_4)\hspace{2mm}\text{for}\hspace{2mm}(s,t)\in G.
\end{equation}
Similarly, if $S$ is a real member, the induced real structure  on $H^0(-2K_S)^{\vee}$ can be seen to be of the form \eqref{0031} in the above  coordinate $(x_0,\cdots,x_4)$.
In particular we obtain (iv).
Going back to a general $S$,
since $\phi:S\ra\mathbf{P}^4$ is a $G$-equivariant birational morphism, 
the image $\phi(S)$ must be an irreducible surface which is the closure of some $G$-orbit.
In view of \eqref{003} we consider a rational map
$(x_0,x_1,x_2,x_3,x_4)\mapsto (x_1x_2,x_3x_4,x_0^2)$
from $\mathbf{P}^4$ to $\mathbf{P}^2$.
Since  $G$ acts trivially on the target space,  every $G$-orbit in $\mathbf P^4$ is contained in some fiber of this rational map.
Namely, any  $G$-orbit is contained in the surface 
\begin{equation}\label{005}
(x_1x_2,x_3x_4,x_0^2)=(a,b,c)
\end{equation}
for some $(a,b,c)\in\mathbf{P}^2$.
Because \eqref{005} always defines a surface of degree 4 in $\mathbf{P}^4$, $\phi(S)$ must coincide with the surface \eqref{005} for some $(a,b,c)$.
Now we show that  $abc\neq 0$.
Since $\phi(S)$ is not contained in any hyperplane in $\mathbf{P}^4$, we have $c\neq 0$. 
If $a=0$ and $bc\neq 0$, then the surface \eqref{005} becomes $\{x_1x_2=0, cx_3x_4=bx_0^2\}$ and this is a reducible quartic surface whose  irreducible components are quadratic cones in $\mathbf{P}^3$ over a smooth conic.
This cannot happen since $\phi(S)$ is irreducible and degree 4.
Similarly, the case $b=0$, $ac\neq 0$ cannot happen.
Moreover,  if $a=b=0$, then the surface \eqref{005} becomes $\{x_1x_2=x_3x_4=0\}$  and this consists of   4 planes.
Hence this case is also impossible.
Thus we obtain $abc\neq 0$.
Letting $c=1$, we obtain the required defining equation \eqref{002} and we get (iii).
\proofend

\vsp
For the proof of Lemma \ref{lemma-011} (i), it is convenient to take a projection from 
a point $(1,0,0,0,0)\in\mathbf P^4$ to a hyperplane $\mathbf{P}^3_{\infty}=\{x_0=0\}$.
(ii) of the lemma can be verified by direct computations using \eqref{002} and \eqref{0031}.
(iii) can also be proved readily because as in the above proof of Lemma \ref{lemma-01} we know the movable part of $|-2K_S|$ explicitly.
Details are left to the reader.

\vsp\noindent
{\bf Proof of (i) of Theorem \ref{thm-01}}.
 By the Riemann-Roch formula, we have, for a twistor space of $n\mathbf{P}^2$, 
$
\chi((-1/2)K_Z)=10-2n.
$
Further, by the positivity of the scalar curvature of the Joyce metric, we have $H^2((-1/2)K_Z)=0$ by Hitchin's vanishing theorem \cite{Hi80}.
By the Serre duality we also have $H^3((-1/2)K_Z)=0$.
Hence  if $n=4$ we have $\dim H^0((-1/2)K_Z)-\dim H^1((-1/2)K_Z)=2$ and it follows $H^1((-1/2)K_Z)=0$.
Next take an irreducible member  $S\in |(-1/2)K_Z|$ and  consider an exact sequence 
\begin{equation}\label{008}
0\lra (-1/2)K_Z\lra -K_Z\lra -2K_S\lra 0.
\end{equation}
From this we obtain an exact sequence 
\begin{equation}\label{013}
0\lra H^0((-1/2)K_Z)\lra H^0(-K_Z)\lra H^0(-2K_S)\lra 0.
\end{equation}
Since we have $\dim H^0(-2K_S)=5$  by Lemma \ref{lemma-01} (i), 
we obtain $\dim H^0(-K_Z)=7$.
Finally, since $2S\in |-K_Z|$ and by the surjectivity of the restriction map $H^0(-K_Z)\ra H^0(-2K_S)$, we have Bs\,$|-K_Z|=$\, Bs\,$|-2K_S|$.
But  the latter  is contained in the cycle $C$ by Lemma \ref{lemma-01} (i).
\proofend

\vsp
\noindent
{\bf Proof of (ii) of Theorem \ref{thm-01}}.
As in the theorem let 
$\Phi:Z\ra\mathbf P^{\vee}H^0(-K_Z)=\mathbf{P}^6$ be the meromorphic map associated to the anticanonical system of $Z$.
The goal is to obtain a diagram \eqref{016} below, which immediately implies the bimeromorphicity of $\Phi$.
To this end we first show that the $G$-action on $H^0(-K_Z)$ induced from  $Z$ is given by 
\begin{equation}\label{014}
(y_1,y_2,y_3,y_4,y_5,y_6,y_7)\longmapsto (s^{-1}y_1,sy_2,t^{-1}y_3,ty_4,y_5,y_6,y_7)\hspace{2mm}\text{for}\hspace{2mm}(s,t)\in G
\end{equation}
for a coordinate $(y_1,\cdots,y_7)$ on $H^0(-K_Z)$.
For this, by \eqref{now} and the $G$-equivariant exact sequence \eqref{013}, it suffices to show that the natural $G$-action on $H^0((-1/2)K_Z)\simeq\mathbf C^2$ is trivial.
Since each member $S\in |(-1/2)K_Z|$ is $G$-invariant, the $G$-action on $H^0((-1/2)K_Z)$ must be of the form $(x,y)\mapsto (s^mt^nx, s^mt^ny)$ for some $m,n\in \mathbf Z$.
To show $m=n=0$, let $V\subset H^0(-K_Z)$ be  the subspace generated by the image of the map 
\begin{equation}\label{0141}
H^0((-1/2)K_Z)\times H^0((-1/2)K_Z)\lra H^0(-K_Z);\hspace{3mm}
(\sigma_1,\sigma_2)\longmapsto \sigma_1\otimes\sigma_2.
\end{equation}
$V$ is clearly 3-dimensional, $G$-invariant and the $G$-action is a scalar multiplication of $s^{2m}t^{2n}$.
By the exact sequence \eqref{013}, the image of $V$ under the restriction map $H^0(-K_Z)\ra H^0(-2K_S)$ is 1-dimensional, on which $G$ still acts by multiplication of $s^{2m}t^{2n}$.
By \eqref{now} this can happen only when $m=n=0$, and thus we obtain \eqref{014}.
At the same time we obtain an important fact that the subspace $V$ is the $G$-fixed part of $H^0(-K_Z)$.

By the inclusion $V\subset H^0(-K_Z)$ we obtain  a projection $$\pi:\mathbf P^{\vee}H^0(-K_Z)\,\lra\, \mathbf P^{\vee}V=\mathbf{P}^2.$$
Let $\mathbf{P}^3_{\infty}\subset \mathbf P^{\vee}H^0(-K_Z)$ be the center of $\pi$.
Namely, $\mathbf{P}^3_{\infty}$ is the projective space of the linear space $\{l\in H^0(-K_Z)^{\vee}\set \,l|_{V}=0\}$.
Fibers of $\pi$ are 4-dimensional linear subspaces in $\mathbf P^{\vee}H^0(-K_Z)$ containing $\mathbf{P}^3_{\infty}$.
For simplicity of notation we write $W_S=H^0(-2K_S)$, which is 5-dimensional.
Since the restriction map $H^0(-K_Z)\ra W_S$ is surjective,
 $\Phi|_S$ can be identified with the morphism $\phi:S\ra \mathbf P^{\vee}W_S$ in Lemma \ref{lemma-01},
 where 
 $ \mathbf P^{\vee}W_S$ can be canonically regarded as a linear subspace of $\mathbf P^{\vee}H^0(-K_Z)$ by the surjection.
 Since the image of the inclusion $H^0((-1/2)K_Z)\ra H^0(-K_Z)$ in \eqref{013} is contained in $V$, 
 $\mathbf P^3_{\infty}$ is contained in $\mathbf P^{\vee}W_S$ for any irreducible $S\in |(-1/2)K_Z|$.
 Hence the image $\pi(\mathbf P^{\vee}W_S)$ is a point for each $S\in |(-1/2)K_Z|$.
 In this way the pencil $|(-1/2)K_Z|$ determines a holomorphic curve $\Lambda$ in $\mathbf P^{\vee}V$.
We  show that $\Lambda$ is a smooth conic.
 To see this, we note that the subspace $\mathbf P^{\vee}W_S$ in $\mathbf P^{\vee}H^0(-K_Z)$ is determined by the kernel of the surjection $H^0(-K_Z)\ra W_S$, and that $\pi (\mathbf P^{\vee}W_S)$  is precisely the point of $\mathbf P^{\vee}V$ determined by   the kernel.
 Concretely, let $\{\sigma_1,\sigma_2\}$ be a basis of $H^0((-1/2)K_Z)$.
Then if $S\in  |(-1/2)K_Z|$ is defined by $u_1\sigma_1+u_2\sigma_2$ for some $(u_1,u_2)\in\mathbf C^2$, the image of the inclusion $H^0((-1/2)K_Z)\ra H^0(-K_Z)$ is the plane $\{(v_1\sigma_1+v_2\sigma_2)\otimes (u_1\sigma_1+u_2\sigma_2)\set (v_1,v_2)\in\mathbf C^2\}$.
The latter can be rewritten as $\{u_1v_1\sigma_1^2+u_2v_2\sigma_2^2+(u_2v_1+u_1v_2)\sigma_1\sigma_2\set (v_1,v_2)\in\mathbf C^2\}$ whose Pl\"ucker coordinate is given by $(u_1u_2,-u_2^2,-u_1^2)$.
It follows that $\Lambda$ is a smooth conic in $\mathbf P^{\vee}V$.
From the construction we have the following commutative diagram of meromorphic maps:
 \begin{equation}\label{016}
 \CD
 Z@>\Phi>> \mathbf P^6\\
 @V\Psi VV @VV\pi V\\
 \Lambda@>\iota>> \mathbf P^2
 \endCD
 \end{equation}
 where $\Psi:Z\ra\Lambda$ is the meromorphic map associated to the pencil $|(-1/2)K_Z|$ and
$\iota:\Lambda\ra \mathbf P^2=\mathbf P^{\vee}V$ is the inclusion as the conic.
$\Lambda$ must be real  since the pencil $|(-1/2)K_Z|$ is  real.
As is already seen, if  $S_{\lambda}:=\Psi^{-1}(\lambda)$ for $\lambda\in \Lambda$ is irreducible, $\Phi|_{S_{\lambda}}$ is just the morphism associated to $|-2K_{S_{\lambda}}|$ and is birational onto its image by Lemma \ref{lemma-01} (ii).
Since general members of $|(-1/2)K|$ are irreducible and $\iota$ is  an embedding, it follows that $\Phi$ is bimeromorphic onto its image.
Thus we obtain (ii) of Theorem \ref{thm-01}.
\proofend

\vsp
\noindent
{\bf Proof of (iii) and (iv) of Theorem \ref{thm-01}.}
Let $(x_1,\cdots,x_7)$ be the coordinate on $H^0(-K_Z)^{\vee}$ which is dual to the coordinate $(y_1,\cdots,y_7)$ for which the natural $G$-action is given by \eqref{014}.
In this coordinate, the induced $G$-action on $H^0(-K_Z)^{\vee}$ is given by
\begin{equation}\label{017}
(x_1,x_2,x_3,x_4,x_5,x_6,x_7)\longmapsto (sx_1,s^{-1}x_2,tx_3,t^{-1}x_4,x_5,x_6,x_7)\hspace{2mm}\text{for}\hspace{2mm}(s,t)\in G.
\end{equation}
It is also readily seen by choosing a real $S$ in \eqref{013} and by \eqref{0031} that the induced real structure on $H^0(-K_Z)^{\vee}$ is given as in \eqref{002}.
In particular we obtain (iv) of the theorem.
Let $V\subset H^0(-K_Z)$, $\mathbf P^3_{\infty}\subset\mathbf P^{\vee}H^0(-K_Z)$, and $\Lambda\subset\mathbf P^2$ be as in the proof of  (ii).
Then $(x_5,x_6,x_7)$ becomes a coordinate on $V^{\vee}$ and we have 
$
\mathbf P^3_{\infty}=\left\{x_5=x_6=x_7=0\right\}.
$
Since $\Lambda\subset\mathbf P^{\vee}V$ is a smooth real conic having real locus, we can suppose that 
\begin{equation}\label{018}
\Lambda=\{(x_5,x_6,x_7)\in\mathbf P^2\set x_5^2=x_6x_7\},
\end{equation}
which is one of the required equations \eqref{001}.
Next we define 
$$
B=\{(x_1,\cdots,x_7)\in \mathbf P^6\set x_5=x_6=x_7=0, \,x_1x_2=x_3x_4=0\}
$$
which is a cycle of 4 lines in $\mathbf P^3_{\infty}$.
$B$ is characterized by the property that it is a set of all non-real $G$-invariant lines in $\mathbf P^3_{\infty}$.

Let $X=\Phi(Z)\subset\mathbf P^6$ be the image of $Z$ and write $X_{\lambda}=\Phi(S_{\lambda})$ for $\lambda\in \Lambda$.
By the diagram  \eqref{016} $X_{\lambda}$ is contained in $\pi^{-1}(\lambda)\simeq\mathbf P^4$ (where we omit the inclusion  $\iota$) which always contains $\mathbf P^3_{\infty}$.
Moreover, $B$ is contained in $\Phi(S_{\lambda})$ for all $\lambda\in\Lambda$, since by Lemmas \ref{lemma-01} and \ref{lemma-011} the 4 components of $B$ must be the images of  $C_2,C_5,\ol{C}_2,\ol{C}_5$, and since $C$ is the base curve of $|(-1/2)K_Z|$.
Thus while $|(-1/2)K_Z|$ is a pencil  whose general members are smooth toric surface having $C$ as the base locus, their image under $\Phi$ is a pencil whose general members are toric surface with 4 ordinary double points, having $B$ as the base locus.

To obtain defining equations of $X$ in $\mathbf P^6$, we blow up $\mathbf P^6$ along $\mathbf P^3_{\infty}$. 
Then the projection $\pi:\mathbf P^6\ra\mathbf P^2$ becomes a morphism and $\mathbf P^6$ becomes the total space of the $\mathbf P^4$-bundle 
$$
p:\mathbf P(\mathscr O(1)^{\oplus 4}\oplus \mathscr O)\lra\mathbf P^2,
$$
whose fibers are naturally identified with fibers of $\pi$.
Define $X_1\subset \mathbf P(\mathscr O(1)^{\oplus 4}$ to be the strict transform of $X$, and let $\nu_1:X_1\ra X$ be the natural surjection and $p_1:X_1\ra\Lambda$ the restriction of $p$ onto $X_1$.
$\nu_1$ is the blowing-up of $X_1$ along the cycle $B$.
Since $\Lambda$ is a conic, $X_1$ is contained in $p^{-1}(\Lambda)\simeq\mathbf P(\mathscr O(2)^{\oplus 4}\oplus \mathscr O)$ over $\Lambda\simeq\mathbf P^1$.
We put
$
\xi_i=x_i/x_7$ for $1\leq i\leq 4$
and $\lambda=x_5/x_7$.
 Then  we can use $\lambda$ as an affine coordinate on $\Lambda\backslash\{(0:1:0)\}\simeq\mathbf C$, and  $(\xi_1,\xi_2,\xi_3,\xi_4)$ as a fiber coordinate of the bundle $\mathscr O(2)^{\oplus 4}\ra\Lambda$ over there.
 Then by \eqref{0011} the natural real structure on the bundle is given by 
 \begin{equation}\label{0181}
 (\xi_1,\xi_2,\xi_3,\xi_4;\lambda)\longmapsto (\ol{\xi}_2,\ol{\xi}_1,\ol{\xi}_4,\ol{\xi}_3;\ol{\lambda}).
 \end{equation}
 Since each $X_{\lambda}$ ($\lambda\in \Lambda$) contains the center $B$, $\nu_1$ is isomorphic on the strict transform of $X_{\lambda}$, which is contained in $p^{-1}(\lambda)$.
If we use $(\xi_1,\xi_2,\xi_3,\xi_4)$ as an affine fiber coordinate on the bundle $\mathbf P(\mathscr O(2)^{\oplus 4}\oplus \mathscr O)$,  by  the equations \eqref{002}, the strict transform of $X_{\lambda}$ in $X_1$ is defined by
\begin{equation}\label{019}
\xi_1\xi_2=f_1(\lambda),\,\,\,\xi_3\xi_4=f_2(\lambda)
\end{equation}
where $f_1(\lambda)$ and $f_2(\lambda)$ are holomorphic functions of $\lambda\in\Lambda$.
 But since \eqref{019} is an equation with values in $\mathscr O(4)$, $f_1$ and $f_2$ must be a quartic polynomial of $\lambda$. 
 Moreover, by Lemma \ref{lemma-01} and \eqref{0181}, $f(\lambda)$ must be real number as long as $\lambda$ is real.
 Let $\Lambda_1\subset\Lambda$ and  $\Lambda_2\subset\Lambda$ be the zero locus of the equations $f_1(\lambda)=0$ and $f_2(\lambda)=0$ respectively and
 we view $\Lambda$ as the conic \eqref{018} in $\mathbf P^2$.
 Since $\Lambda$ is real and $\Lambda_i$ ($i=1,2$) are real subsets, there exist real conics  $\mathscr C_i$ ($i=1,2$) such that $\mathscr C_i\cap\Lambda=\Lambda_i$ hold scheme theoretically.
 (The choices of $\mathscr C_1$ and $\mathscr C_2$ are not unique.)
  Let $Q_1(x_5,x_6,x_7)$ and $Q_2(x_5,x_6,x_7)$ be defining quadratic polynomials of $\mathscr C_1$ and $\mathscr C_2$ respectively, which  necessarily have real coefficients by the choice. 
 Then by the definition of $\mathscr C_i$, we have, for $i=1,2$, 
 \begin{equation}\label{020}
f_i\left(\frac{x_5}{x_7}\right)=c_i\,Q_i\left(\frac{x_5}{x_7},\left(\frac{x_5}{x_7}\right)^2,1\right)
 \end{equation}
 for some non-zero constants $c_1$ and $c_2$.
Therefore by substituting $\lambda=x_5/x_7$ to \eqref{019} and then multiplying them by $x_7^2$ and using $x_5^2=x_6x_7$, we obtain the remaining two equations of \eqref{001}.
\proofend

\vsp
The structure of the bimeromorphic map $\Phi$ can be understood more readily if we blowup $Z$. 
Namely if $\mu_1:\tilde Z\ra Z$ denotes the blowing-up along the cycle $C$, and if $\nu_1:X_1\ra X$ still denotes the blowing-up along $B$ as in the proof of (iii), then  by \eqref{016} we obtain a commutative diagram:
\begin{equation}\label{100}
 \CD
 \tilde Z@>\tilde{ \Phi}>> X_1\\
 @VVV @VVV\\
 \Lambda@= \Lambda,
 \endCD
 \end{equation}
 where vertical arrows are morphisms but $\tilde\Phi$ (which is induced from $\Phi$) still has  indeterminacy.
Note that both $X_1$ and $\tilde Z$ have singularities. 
As in the proof of (iii), $X_1$ is naturally a subvariety in the $\mathbf P^4$-bundle $\mathbf P(\mathscr O(2)^{\oplus 4}\oplus \mathscr O)$ over $\Lambda\simeq\mathbf P^1$ and its defining equations are given by \eqref{019}.

\subsection{More restrictions on the defining equations}
Next we will obtain more necessary conditions for $Q_1$ and $Q_2$ in Theorem \ref{thm-01}, by taking reducible members of the pencil $|(-1/2)K_Z|$ into account.
To this end, as in the proof of (ii) and (iii) of Theorem \ref{thm-01} we consider a projective plane $\mathbf P^{\vee}V=\mathbf P^2$ having $(x_5:x_6:x_7)$ as a homogeneous coordinate 
and write $\mathscr C_i$ ($i=1,2$) for the conics in the plane defined by $Q_i(x_5,x_6,x_7)=0$ respectively, and denote $\Lambda_i=\Lambda\cap\mathscr C_i$ $(i=1,2)$ as before.
Of course, these 3 conics are respectively real with respect to the real structure $(x_5,x_6,x_7)\mapsto (\ol{x}_5,\ol{x}_6,\ol{x}_7)$.
Moreover, if $\Lambda^{\sigma}$ denotes the set of real points in $\Lambda$, then $\Lambda^{\sigma}\simeq S^1$.

\begin{figure}
\includegraphics{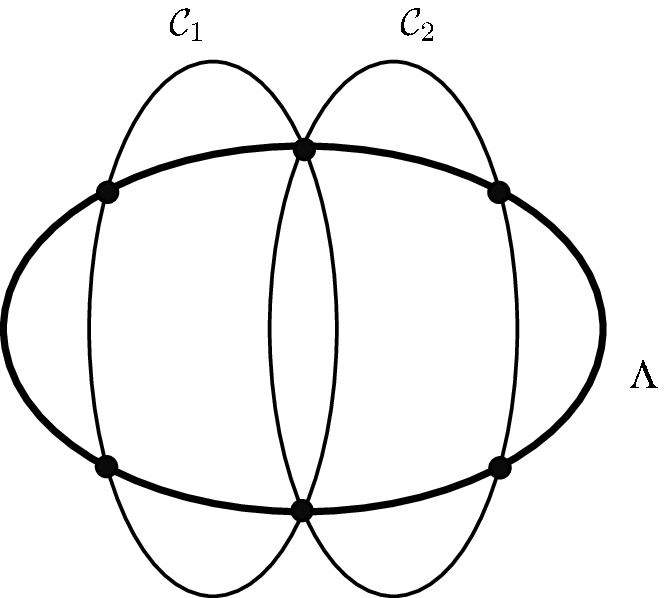}
\caption{The 3 conics determining the projective model of the twistor space}
\label{fig-3conics}
\end{figure}

\begin{thm}
\label{thm-02} For the intersections of these conics we have the following.
(i) $\Lambda_1\cap\Lambda_2$ consists of 2 real points.
(ii) Both $\Lambda_1$ and $\Lambda_2$ consist of 4 real points, so that $\Lambda_1\cup\Lambda_2$ consists of 6 real points.
(iii) The 2 points $\Lambda_1\backslash(\Lambda_1\cap\Lambda_2)$ belong to one of the 2 connected components of  $\Lambda^{\sigma}\backslash(\Lambda_1\cap\Lambda_2)$, and 
the 2 points $\Lambda_2\backslash(\Lambda_1\cap\Lambda_2)$ belong to another component of $\Lambda^{\sigma}\backslash(\Lambda_1\cap\Lambda_2)$.
\end{thm}
This says that the 3 conics intersect as illustrated in Figure \ref{fig-3conics}.

\begin{figure}
\includegraphics{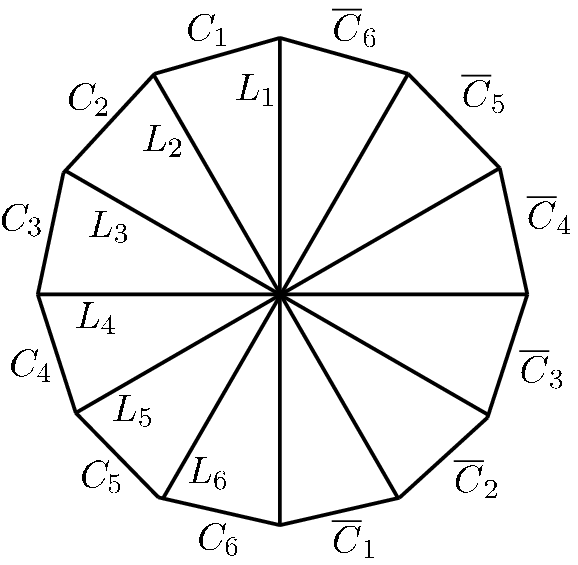}
\caption{The numbering for $G$-invariant twistor lines}
\label{fig-6tls}
\end{figure}

In the rest of this subsection we will prove   Theorem \ref{thm-02}.
For this we need to know the images of the reducible members of the pencil $|(-1/2)K_Z|$ under the bimeromorphic map $\Phi$, basically because the zeros of $f_1$ and $f_2$ in the equations \eqref{739} correspond to the reducible members of the pencil.
So we give a numbering for $G$-invariant twistor lines as  illustrated in Figure \ref{fig-6tls}.
Then recalling that each $L_i$ ($1\le i\le 6$) is of the form $S_i^+\cap S_i^-$ (Proposition \ref{prop-01} (v)), we obtain a natural numbering for the reducible members of the pencil.
Moreover, though there is no natural way to make a distinction between $S_i^+$ and $S_i^-$, we define $S_i^+$ $(1\le i\le 6)$ to be the component which contains $C_6$, so that $S_i^-$ is the component containing $\ol{C}_6$.

\begin{lemma}\label{lemma-10}
The following divisors are all $G$-invariant members of $|-K_Z|$:
(a) $S+S'$, where $S,\,S'\in|(-1/2)K_Z|$,
(b) $S_1^++S_2^++S_3^++S_4^-$, $S_1^-+S_2^-+S_3^-+S_4^+$, $S_1^++S_4^++S_5^++S_6^+$ and $S_1^-+S_4^-+S_5^-+S_6^-$ .
\end{lemma}
\noindent Proof.
As seen in \eqref{014}, the natural $G$-action on $H^0(-K_Z)=\mathbf C^7$  is of the form
$$(y_1,y_2,y_3,y_4,y_5,y_6,y_7)\longmapsto (s^{-1}y_1,sy_2,t^{-1}y_3,ty_4,y_5,y_6,y_7)
$$
 for a coordinate on $H^0(-K_Z)$.
Hence if $l$ is a $G$-invariant 1-dimensional subspace of $H^0(-K_Z)$, then one of the following holds:
$(a')$\, $l$ is contained in $V=\{y_1=y_2=y_3=y_4=0\}$, $(b')$ $l$ is generated by $(1,0,0,0,0,0), (0,1,0,0,0,0,0), (0,0,1,0,0,0,0)$ or $(0,0,0,1,0,0,0)$.
As in the proof of (ii) of Theorem \ref{thm-01}, $V$ is generated by the image of the map \eqref{0141}.
Hence if $l$ is in $(a')$, then $l$ must be of the form $(a)$ in the lemma.
It remains to see that the 4 divisors in $(b)$ in the lemma are actually members of $|-K_Z|$.
This can be verified by explicitly specifying the cohomology classes of $S_i^{\pm}$ in $H^2(Z,\mathbf Z)$ for all $1\le i\le 6$ and checking that the sums of the 4 cohomology classes belong to $c_1(Z)$.
We omit the detail.
\proofend

\vsp
Next for each $1\le i\le 6$, we define $W^{\pm}_i$ to be the images of the restriction $H^0(-K_Z)\ra H^0(-K_Z|S_i^{\pm})$, so that we have an exact sequence
\begin{equation}
\label{101}
0\,\lra H^0(-K_Z-S_{i}^{\pm})\,\lra H^0(-K_Z)\,\lra\,W_i^{\pm}\,\lra 0.
\end{equation}
Further we write $V_i^{\pm}\,(\subset H^0(-K_Z))$ for the images of the injection in \eqref{101}.
If $|W_i^{\pm}|$ denote the linear systems on $S_i^{\pm}$ associated to $W_i^{\pm}$, the restriction $\Phi|{S_i^{\pm}}$ is precisely the associated rational map $S_i^{\pm}\ra \mathbf P^{\vee}W_i^{\pm}$ of $|W_i^{\pm}|$, where $\mathbf P^{\vee}W_i^{\pm}$ are naturally   linear subspaces of $\mathbf P^{\vee}H^0(-K_Z)=\mathbf P^6$ by the surjections in \eqref{101}.
Then the following result concerning how these linear subspaces intersect is essential in our proof of Theorem \ref{thm-02}.
Also, the calculations in the proof provide a partial reason for the complicated construction of the twistor spaces given in the next section.

\begin{lemma}\label{lemma-11}
We have the following.
(i) $\dim W_i^{\pm}=3$ for $i=1,4$ and $\dim W_i^{\pm}=4$ for $i=2,3,5,6$.
(ii) $\dim (\mathbf P^{\vee}W_i^+\cap \mathbf P^{\vee}W_i^-)=0$ for $i=1,4$
and $\dim (\mathbf P^{\vee}W_i^+\cap \mathbf P^{\vee}W_i^-)=2$ for $i=2,3,5,6$.
(iii) $\dim(\mathbf P^{\vee}W_1^{\pm}\cap\mathbf P^{\vee}W_4^{\pm})=0$, where ${\pm}$ take every combinations.
(iv)  $\dim(\mathbf P^{\vee}W_2^+\cap\mathbf P^{\vee}W_3^+)=\dim(\mathbf P^{\vee}W_5^+\cap\mathbf P^{\vee}V_6^+)=2$.
\end{lemma}

 \noindent Proof.
 For (i), we have only to prove $\dim W_i^+=3$ thanks to the real structure.
 Since $\dim H^0(-K_Z)=7$, we have $\dim W_i^{+}=7-\dim H^0(-K_Z-S^{+}_i)$. 
 To calculate the latter cohomology group, we consider the exact sequence
 \begin{equation}\label{0191}
 0\,\lra\,-K_Z-S^{+}_i-{S}^{-}_i\,\lra\, -K_Z-S_i^{+}\,\lra\, -K_Z-S^{+}_i|_{{S}_i^{-}}\,\lra\,0.
 \end{equation}
Then since $S_i^{+}+S^{-}_i\in |(-1/2)K_Z|$ we have $-K_Z-S_i^{+}-S_i^{-}\simeq (-1/2)K_Z$.
 On the other hand we have $H^0((-1/2)K_Z)\simeq\mathbf C^2$ and  $H^1((-1/2)K_Z)=0$
 as seen in the proof of Theorem \ref{thm-01} (i). 
 Hence by the cohomology exact sequence of \eqref{0191} we obtain 
 $\dim H^0(-K_Z-S_i^{+})=2+\dim H^0(-K_Z-S_i^{+}|_{S_i^{-}}).$
 Hence we obtain 
 \begin{equation}\label{021}
 \dim W_i^{+}=5-\dim H^0\left({S}_i^{-},(-K_Z-S_i^{+})|_{{S}_i^{-}}\right).
 \end{equation} 
Next we calculate the right hand side of \eqref{021}.
By adjunction we have $K_{S^{-}_i}\simeq K_Z|_{S^{-}_i}+N_{S^{-}_i/Z}.$
On the other hand, we have $(-1/2)K_Z|_{S_i^{-}}\simeq S^{+}_i+S^{-}_i|_{S_i^{-}}\simeq N_{S_i^{-}/Z}+L_i$.
From these we readily obtain
$(-1/2)K_Z|_{S_i^{-}}\simeq -K_{S^{-}_i}-L_i.$
Then noting $S_i^{+}|_{S_i^{-}}\simeq L_i$, we obtain
$
\left(-K_Z-S_i^{+}\right)|_{S_i^{-}}\simeq -2K_{S_i^{+}}-3L_i .
$
Now since $S^{+}_i$ is a toric surface whose structure can be explicitly described (cf an explanation after Prop.\,\,\ref{prop-01}),
we readily obtain an effective (and torus invariant) member of $|-2K_{S_i^{+}}-3L_i|$.
Once a member is obtained it is not difficult to compute the dimension of the associated linear system, for each $S^{+}_i$, $1\le i\le 6$.
We omit the detail since it is a succession of standard computations.
As a result we obtain that for $i=1,4$ we have $\dim H^0(-K_Z-S^{+}_i|_{S_i^{-}})=2$ and the remaining 4 cases we have $\dim H^0(-K_Z-S_i^{+}|_{S_i^{-}})=1$.
Hence by \eqref{021} we obtain (i) of the lemma.

For (ii) we note that $\mathbf P^{\vee}W_i^{\pm}$ is a linear subspace associated to the vector space $\{l\in H^0(-K_Z)^{\vee}\set\,\, l|_{V_i^{\pm}}\equiv 0\}$.
From this it follows that 
\begin{equation}\label{0231}
\dim (\mathbf P^{\vee} W_i^+\cap\mathbf P^{\vee} W_i^-)=
\dim H^0(-K_Z)-\dim(V_i^++V_i^-)-1
\end{equation}
where $V_i^++V_i^-$ is the usual sum as vector subspaces in $H^0(-K_Z)$.
On the other hand by the natural isomorphisms $V_i^+\cap V_i^-\simeq H^0(-K_Z-S_i^+-S_i^-)\simeq H^0((-1/2)K_Z)\simeq\mathbf C^2$ we have $\dim (V_i^+\cap V_i^-)=\dim V_i^++\dim V_i^--2$.
Hence by (i) and \eqref{0231} we obtain (ii) of the lemma.

For (iii), by reality it is enough to show that $\dim (\mathbf P^{\vee}W_1^+\cap\mathbf P^{\vee}W_4^+)=\dim (\mathbf P^{\vee}W_1^+\cap\mathbf P^{\vee}W_4^-)=0$.
For this using the same argument as in (ii) it suffices to verify that $\dim (V_1^+\cap V_4^+)=\dim (V_1^+\cap V_4^-)=2$, or equivalently $\dim H^0(-K_Z-S_1^+-S_4^+)=\dim H^0(-K_Z-S_1^+-S_4^-)=2$.
To see the latter, first note that $S_1^-+S_4^-$ is a member of $|-K_Z-S_1^{+}-S_4^{+}|$.
Moreover, by Lemma \ref{lemma-10}, $S_5^++S_6^+$ is also a member of $|-K_Z-S_1^{+}-S_4^{+}|$.
From these we obtain $\dim H^0(-K_Z-S_1^+-S_4^+)\ge 2$.
Moreover, if $\dim H^0(-K_Z-S_1^+-S_4^+)\ge 3$, there must be another $G$-invariant member of $|-K_Z-S_1^{+}-S_4^{+}|$.
Such a member must be a sum of two degree 1 divisors, since all $G$-invariant divisors in $Z$ are members of $|(-1/2)K_Z|$ or $S_{i}^{\pm}$, $1\le i\le 6$, and since $-K_Z-S_1^{+}-S_4^{+}\not\simeq(-1/2)K_Z$.
But again by Lemma \ref{lemma-10} there does not exist  such a member other than $S_1^++S_4^+$ and $S_5^++S_6^+$.
Thus we obtain $\dim H^0(-K_Z-S_1^+-S_4^+)= 2$.
By the same argument it follows that $S_1^-+S_4^+$ and $S_2^++S_3^+$ 
are all $G$-invariant members of $|-K_Z-S_1^+-S_4^-|$ and we also obtain  $\dim H^0(-K_Z-S_1^+-S_4^-)=2$.  
Hence we obtain (iii) of the lemma.

To show $\dim (\mathbf P^{\vee}W_2^+\cap\mathbf P^{\vee} W_3^+)=2$, by the same reason as in (ii) and (iii), it suffices to show that $\dim (V_2^+\cap V_3^+)=2$ or equivalently $\dim H^0(-K_Z-S_2^+-S_3^+)=2$.
It is clear that $S_2^-+S_3^-$ is a member of $|-K_Z-S_2^+-S_3^+|$  and by Lemma \ref{lemma-10}, $S_1^++S_4^-$ is also a member of $|-K_Z-S_2^+-S_3^+|$.
Hence $\dim H^0(-K_Z-S_2^+-S_3^+)\ge 2$.
To obtain the equality, we can argue in the same way as in (iii).
$\dim (\mathbf P^{\vee}W_5^+\cap\mathbf P^{\vee} W_6^+)=2$ can be also verified exactly in the same way.
Thus we obtain (iv) of the lemma.
\proofend

 \vsp
\noindent
{\bf Proof of Theorem \ref{thm-02}.}		
We freely use notations in the proof of Theorem \ref{thm-01}.
For $\lambda\in\Lambda$, $\Phi(S_{\lambda})$ is contained in a surface in $\pi^{-1}(\lambda)\simeq\mathbf P^4$ (cf. the diagram \eqref{016}) defined by 
\begin{equation}\label{383}
\xi_1\xi_2=f_1(\lambda),\,\,\,\xi_3\xi_4=f_2(\lambda).
\end{equation}
Let $\lambda_i\in\Lambda$, $1\le i\le 6$, be the points such that $S_i=S_{\lambda_i}$ holds.
If $\lambda\neq \lambda_i$ for $1\le i\le 6$, namely if $S_{\lambda}$ is irreducible, then as showed in the proof Theorem \ref{thm-01},  $\Phi|_{S_{\lambda}}$ coincides with the map associated to $|-2K|$ and hence by Lemma \ref{lemma-01}, $\Phi(S_{\lambda})$ is the surface \eqref{002} with $ab\neq 0$. 
Therefore we have $f_1(\lambda)f_2(\lambda)\neq 0$ for $\lambda\neq\lambda_i$.
Suppose $i=1$ or $i=4$.
 $\Phi(S^{\pm}_i)$ is contained in  $\mathbf P^{\vee}W^{\pm}_i\subset\mathbf P^{\vee}H^0(-K_Z)$ and  $\mathbf P^{\vee}W^{\pm}_i$ is a plane by Lemma \ref{lemma-11} (i), and it is clearly $G$-invariant.
Hence the surface \eqref{383} must contain a $G$-invariant plane.
By \eqref{017}, the $G$-action on $\pi^{-1}(\lambda)$ is of the form $(\xi_1,\xi_2,\xi_3,\xi_4)\mapsto (s\xi_1,s^{-1}\xi_2,t\xi_3,t^{-1}\xi_4)$.
Hence all $G$-invariant planes in $\pi^{-1}(\lambda)$ are $\{\xi_1=\xi_2=0\}$, $\{\xi_1=\xi_3=0\}$,
$\{\xi_1=\xi_4=0\}$, $\{\xi_2=\xi_3=0\}$, $\{\xi_2=\xi_4=0\}$ and $\{\xi_3=\xi_4=0\}$.
Since the real structure on $\pi^{-1}(\lambda)$ ($\lambda\in\mathbf R$) is given by $(\xi_1,\xi_2,\xi_3,\xi_4)\mapsto (\ol{\xi}_2,\ol{\xi}_1,\ol{\xi}_4,\ol{\xi}_3)$ as in \eqref{0181},
among these 6 planes, the first and the last ones are real.
If $\mathbf P^{\vee}W^+_i$  coincides with one of these 2 real planes, then by reality $\mathbf P^{\vee}W^{-}_i$ must be the same real plane.
This contradicts Lemma \ref{lemma-11} (ii).
Hence $\{\mathbf P^{\vee}W^+_i,\mathbf P^{\vee}W^-_i\}$ is a conjugate pair of planes which are among the remaining 4 planes.
This immediately implies $f_1(\lambda_i)=f_2(\lambda_i)=0$ for $i=1,4$.
If $i=2,3,5,6$, $\mathbf P^{\vee} W_i^{\pm}=\mathbf P^3$ by Lemma \ref{lemma-11} (i) and
hence the surface \eqref{383} must contain a $G$-invariant hyperplane in $\mathbf P^4$.
This implies that $f_1(\lambda_i)f_2(\lambda_i)=0$.
Then since $\deg f_1=\deg f_2=4$ and we already know that $f_1(\lambda_i)=f_2(\lambda_i)=0$ for $i=1,4$, precisely one of $f_1(\lambda_i)=0$ and $f_2(\lambda_i)=0$ occurs for each $i=2,3,5,6$.
This implies that both $f_1=0$ and $f_2=0$ have no double roots and precisely two of the four roots coincide.
Thus we obtain (i) and (ii) of the theorem.

Finally to prove (iii) it is enough to  show that $f_1(\lambda_2)=0$ iff $f_1(\lambda_3)=0$, and $f_1(\lambda_5)=0$ iff $f_1(\lambda_6)=0$.
If $f_1(\lambda_2)=0$, then by \eqref{383}, $\Phi(S_2^+)$ and $\Phi(S_2^-)$ are contained $\xi_1\xi_2=0$.
Hence either $\mathbf P^{\vee}W_2^+=\{\xi_1=0\}$ or $\mathbf P^{\vee}W_2^+=\{\xi_2=0\}$ holds.
Suppose the former holds. 
Then by Lemma \ref{lemma-11} (iv), $\mathbf P^{\vee}W_3^+=\{\xi_1=0\}$ must hold.
In the same way, if the latter holds, $\mathbf P^{\vee}W_3^+=\{\xi_2=0\}$ must hold.
Thus we always have $f_1(\lambda_3)=0$.
In the same way, if $f_1(\lambda_3)=0$, then $f_1(\lambda_2)=0$.
By the same argument we have $f_1(\lambda_5)=0$ iff $f_1(\lambda_6)=0$.
Thus $\Lambda_1=\{\lambda_1,\lambda_4,\lambda_2,\lambda_3\}$ and  $\Lambda_2=\{\lambda_1,\lambda_4,\lambda_5,\lambda_6\}$ hold simultaneously or otherwise
$\Lambda_1=\{\lambda_1,\lambda_4,\lambda_5,\lambda_6\}$ and  $\Lambda_2=\{\lambda_1,\lambda_4,\lambda_2,\lambda_3\}$ hold simultaneously.
This directly implies the conclusion of (iii) of Theorem \ref{thm-02}, as required.
\proofend

\vsp
Finally in this section we prove by using Lemma \ref{lemma-11} some detailed result about the bimeromorphic map $\Phi$ needed in the final part of the next section.
For $\lambda\in \Lambda$ we denote by $T_{\lambda}\subset S_{\lambda}$ for the union of all 2-dimensional $G$-orbits in $S_{\lambda}$.
If $S_{\lambda}$ is irreducible, $T_{\lambda}$ is isomorphic to $G$.
If $S_{\lambda}$ is reducible, $T_{\lambda}$ is isomorphic to a (disjoint) union of two copies of $G$. 
Define $U:=\cup_{\lambda\in\Lambda}T_{\lambda}$, the union of all 2-dimensional orbits in $Z$.
In other words, $U$ is the complement of all lower-dimensional $G$-orbits, where the latter consists of 
$G$-orbits contained in $C$ and the $G$-invariant twistor lines $L_i$, $1\le i\le 6$.
Thus $\Phi$ is well defined on $U$ by Theorem \ref{thm-01}.

\begin{lemma}\label{lemma-32}
Let $Z$, $\Phi:Z\ra\mathbf P^6$, $\mathbf P^3_{\infty}\subset \mathbf P^6$, $C\subset Z$, $X=\Phi(Z)\subset\mathbf P^6$ and $B\subset X$ be as before, and $U$ the union of all 2-dimensional $G$-orbits in $Z$ as above.
Then we have the following.
(i) On $U$, $\Phi$ is biholomorphic onto its image and $\Phi(U)$ is contained in $X\backslash B$.
(ii) $\Phi(S_1^+)\cap\mathbf P^3_{\infty}$, $\Phi(S_1^-)\cap\mathbf P^3_{\infty}$, $\Phi(S_4^+)\cap\mathbf P^3_{\infty}$ and $\Phi(S_4^-)\cap\mathbf P^3_{\infty}$ are mutually different non-real $G$-invariant lines in $\mathbf P^3_{\infty}$.
(Hence we have $B=(\Phi(S_1^+)\cup \Phi(S_1^-)\cup \Phi(S_4^+)\cup \Phi(S_4^-))\cap \mathbf P^3_{\infty}$.)
\end{lemma}

We note that the restriction $Z\backslash C\ra X\backslash B$ of $\Phi$ is never surjective, since the meromorphic map  $X\ra\Lambda$ (or the morphism $X_1\ra\Lambda$ in the diagram \eqref{100}) has reducible fibers consisting of 4 planes over $\lambda_i$ for $i=1,4$.
In brief, $X$ has `redundant' divisors which do not exist in the twistor space $Z$.
This is an important and most difficult point in the construction of the twistor spaces given in the next section (especially in Step 5).

\vsp
\noindent Proof of Lemma \ref{lemma-32}.
By Lemma \ref{lemma-011} (iii),
if $\lambda\neq \lambda_i$ for $1\leq i\le 6$, $\Phi|_{S_{\lambda}}$ is a birational morphism which precisely contracts 8 components of $C$.
Thus $\Phi$ gives a biholomorphic map from $S_{\lambda}\backslash C=T_{\lambda}$ to $X_{\lambda}\backslash B$, where $X_{\lambda}=\Phi(S_{\lambda})$.
We show that on the reducible members $S_{i}=S_i^++{S}_i^-$ $(1\le i\le 6)$, $\Phi$ also gives a birational map which is biholomorphic on the two 2-dimensional $G$-orbits.
For $i=1,4$, $\Phi(S_i^{\pm})$ is contained in $\mathbf P^{\vee}W_i^{\pm}=\mathbf P^2$  and the latter must be one of the non-real planes $\{\xi_1=\xi_3=0\}$, $\{\xi_1=\xi_4=0\}$, $\{\xi_2=\xi_3=0\}$ and $\{\xi_2=\xi_4=0\}$ as showed in the proof of (i) and  (ii) of Theorem \ref{thm-02}.
On these 4 planes all $d$-dimensional $G$-orbits with $d=0,1$  are contained in a line.
Therefore by non-degeneracy of the image, it follows $\Phi(S_i^{\pm})=\mathbf P^{\vee}W_i^{\pm}$.
Hence 2-dimensional orbits in $S_i^{\pm}$ must be mapped to 2-dimensional orbits.
Moreover, these must be isomorphic, since $G$ are acting effectively on the planes.
Hence $\Phi$ is isomorphic on $T_{\lambda_i}$ for $i=1,4$.
By similar (and simpler)  argument, it follows that  for $i=2,3,5,6$, $\Phi$ still gives birational maps from $S_i^{\pm}$ onto a quadratic cone in $\mathbf P^{\vee}W_i^{\pm}=\mathbf P^3$ that are biholomorphic on 2-dimensional $G$-orbits.
Thus we have seen that for any $\lambda\in \Lambda$, $\Phi$ is biholomorphic on $T_{\lambda}$.
Then recalling that the pencil $|(-1/2)K_Z|$ has $C$ as the base locus and that they are mapped to a pencil $\{X_{\lambda}\}_{\lambda\in\Lambda}$ which has $B$ as the base locus, it follows that $\Phi$ is biholomorphic on $\cup_{\lambda\in\Lambda}T_{\lambda}=U$ and that $\Phi(U)$ is contained in $X\backslash B$.
Thus we obtain (i).
 
 For (ii), we first note that when $i=1,4$, the set $\{\Phi(S_i^+)\cap\mathbf P^3_{\infty}, \Phi(S_i^-)\cap\mathbf P^3_{\infty}\}$ is  a conjugate pair of $G$-invariant lines in $\mathbf P^3_{\infty}$.
 In fact, if $\Phi(S_i^+)=\mathbf P^{\vee}W_i^+=\mathbf P^2$ were contained in $\mathbf P^3_{\infty}$, then $\Phi(S_i^-)=\mathbf P^{\vee}W_i^-=\mathbf P^2$ would also be contained in $\mathbf P^3_{\infty}$ by reality and therefore $\mathbf P^{\vee}W_i^+\cap\mathbf P^{\vee}W_i^-$ is at least 1-dimensional.
 This contradicts Lemma \ref{lemma-11} (ii).
  Hence both $\Phi(S_i^+)\cap\mathbf P^3_{\infty}$ and  $\Phi(S_i^-)\cap\mathbf P^3_{\infty}$ must be lines, which is obviously a conjugate pair. 
 $G$-invariance of these lines is also clear.
On $\mathbf P^3_{\infty}$ there are precisely 2 pairs of $G$-invariant real lines.
It remains to show that 
 $\{\Phi(S_1^+)\cap\mathbf P^3_{\infty},\Phi(S_1^-)\cap\mathbf P^3_{\infty}\}\cap\{\Phi(S_4^+)\cap\mathbf P^3_{\infty},\Phi(S_4^-)\cap\mathbf P^3_{\infty}\}=\emptyset$.
 If this is not true, then $\{\Phi(S_1^+)\cap\mathbf P^3_{\infty},\Phi(S_1^-)\cap\mathbf P^3_{\infty}\}=\{\Phi(S_4^+)\cap\mathbf P^3_{\infty},\Phi(S_4^-)\cap\mathbf P^3_{\infty}\}$ holds.
 This contradicts  Lemma \ref{lemma-11}, (iii), as required.
 \proofend
 
\section{An explicit construction of the twistor spaces}

In this section we will provide an explicit procedure how to obtain the twistor spaces of Joyce metrics on $4\mathbf P^2$ (whose $K$-action is of Type I) from the projective models  obtained in the previous section.

We first recall continuous parameters involved in the Joyce's construction of his metrics.
A Joyce metric on $n\mathbf P^2$ with a prescribed $K$-action is uniquely determined once the set of $(n+2)$ points $\{y_1<y_2<\cdots<y_{n+2}\}$ on $\mathbf R\cup\{\infty\}=\mathbf{RP}^1$ is given.
Moreover, the metrics coming from different sets of $(n+2)$ points are mutually isomorphic iff they are equivalent under the natural action of $GL(2,\mathbf R)$.
In \cite{F00} the equivalence class of $\{y_1<\cdots<y_{n+2}\}$ is called the  {\em{conformal invariant}} of a Joyce metric.


On the other hand, as in Proposition \ref{prop-01}, if $Z$ is the associated twistor space of a Joyce metric on $4\mathbf{P}^2$,  there exist precisely 6 reducible members $S_1,\cdots, S_6$ of $|(-1/2)K_Z|$. 
Correspondingly we have 6 real points $\lambda_1,\cdots,\lambda_6$ in $\Lambda=\mathbf P^1$ arranged cyclically by our choice of the reducible members  specified in Proposition \ref{prop-01}.
Note that once we give a numbering for $G$-invariant twistor lines as in Figure \ref{fig-6tls},  the numbering for the 6 points in $\Lambda$ is automatically determined, and we adopt this numbering also in this section.
In   \cite[(6.13)]{F00}, an equivalence class of a circular sequence $\{\lambda_1,\cdots,\lambda_6\}$ is called the {\em{twistorial invariant}} of a Joyce metric.
By Theorem 9.1 in \cite{F00}, if $Z$ is the twistor space of a Joyce metric whose conformal invariant is $\{y_1<\cdots<y_6\}$, then  we can suppose that $\lambda_i=y_i$ for $1\le i\le 6$.

Theorems \ref{thm-01} and \ref{thm-02} imply that $Z$ is bimeromorphic to the complete intersection $X$ defined by
\begin{equation}
\label{732}
x_1x_2=Q_1(x_5,x_6,x_7),\,\,
x_3x_4=Q_2(x_5,x_6,x_7),\,\,
x_5^2=x_6x_7,
\end{equation}
where the conics $\mathscr C_i=\{Q_i=0\}$ $(i=1,2)$ intersect $\Lambda=\{x_5^2=x_6x_7\}$ in such a way as in Theorem \ref{thm-02} (see Figure \ref{fig-3conics}).
If we blow up $Z$ along the cycle $C$, then correspondingly $X$ is blown-up along $B=\{x_5=x_6=x_7= x_1x_2=x_3x_4=0\}$ and the strict transform $X_1$ of $X$ was defined by
\begin{equation}\label{739}
\xi_1\xi_2=f_1(\lambda),\,\,\xi_3\xi_4=f_2(\lambda),
\end{equation}
where $\xi_i$ $(1\le i\le 4)$ are fiber coordinates of $\mathscr O(2)$, $f_1$ and $f_2$ are quartic polynomials of $\lambda$, and \eqref{739} is an equation in $\mathbf P(\mathscr O(2)^{\oplus 4}\oplus \mathscr O)\ra\Lambda=\mathbf P^1$.
As in \eqref{017} and \eqref{0181},
the $G$-action and the real structure are respectively given by
$$
(\xi_1,\xi_2,\xi_3,\xi_4;\lambda)\longmapsto (s\xi_1,s^{-1}\xi_2,t\xi_3,t^{-1}\xi_4;\lambda)\,\,\,{\text{for}}\,\,\, (s,t)\in G
$$
and 
$$
(\xi_1,\xi_2,\xi_3,\xi_4;\lambda)\longmapsto (\ol{\xi}_2,\ol{\xi}_1,\ol{\xi}_4,\ol{\xi}_3;\ol{\lambda}).
$$
By virtue of Theorem \ref{thm-02} and Lemma \ref{lemma-011} (ii), we can suppose that 
\begin{align}f_1(\lambda)&=-(\lambda-\lambda_1)(\lambda-\lambda_4)(\lambda-\lambda_2)(\lambda-\lambda_3),\\
f_2(\lambda)&=(\lambda-\lambda_1)(\lambda-\lambda_4)(\lambda-\lambda_5)(\lambda-\lambda_6).
\end{align}
As before let $p_1:X_1\ra \Lambda$ be the projection.
Then the following (a), (b) and (c) are easy to see. (Some of them are already seen in the previous section.)

\begin{figure}
\includegraphics{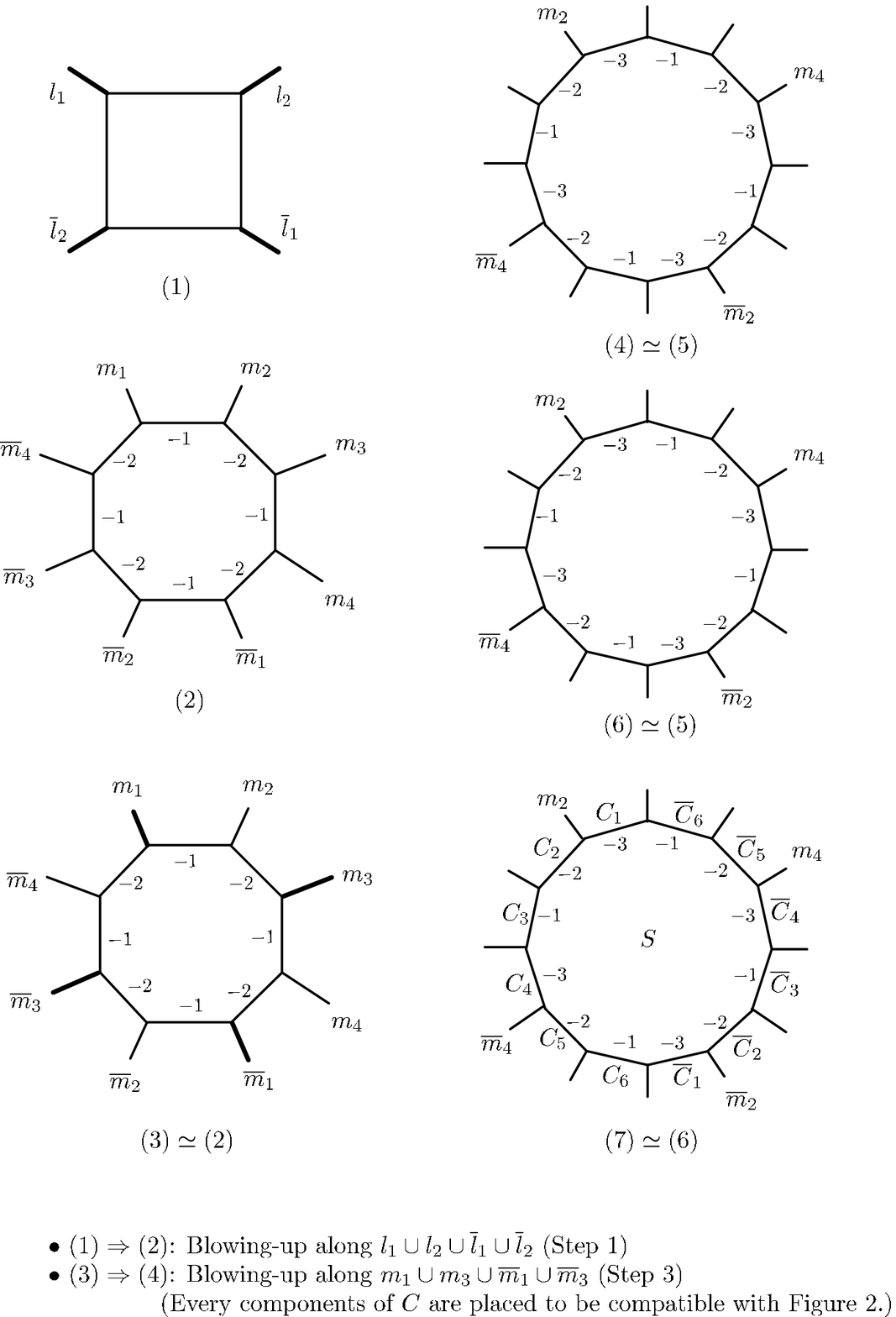}
\caption{The transformations for irreducible fibers}
\label{fig-bim0}
\end{figure}

\vsp\noindent
(a)
Over $\Lambda\backslash\{\lambda_i\set1\le i\le 6\}$, $p_1$ is a fiber bundle whose fibers are irreducible toric surfaces with 4 ordinary double points (Lemma \ref{lemma-011}).
These singularities form 4 global sections of $p_1$, which we call  $l_1,l_2,\ol{l}_1,\ol{l}_2$ (Figure \ref{fig-bim0}\,(1)).
These  are 1-dimensional singular locus of $X_1$.

\begin{figure}
\includegraphics{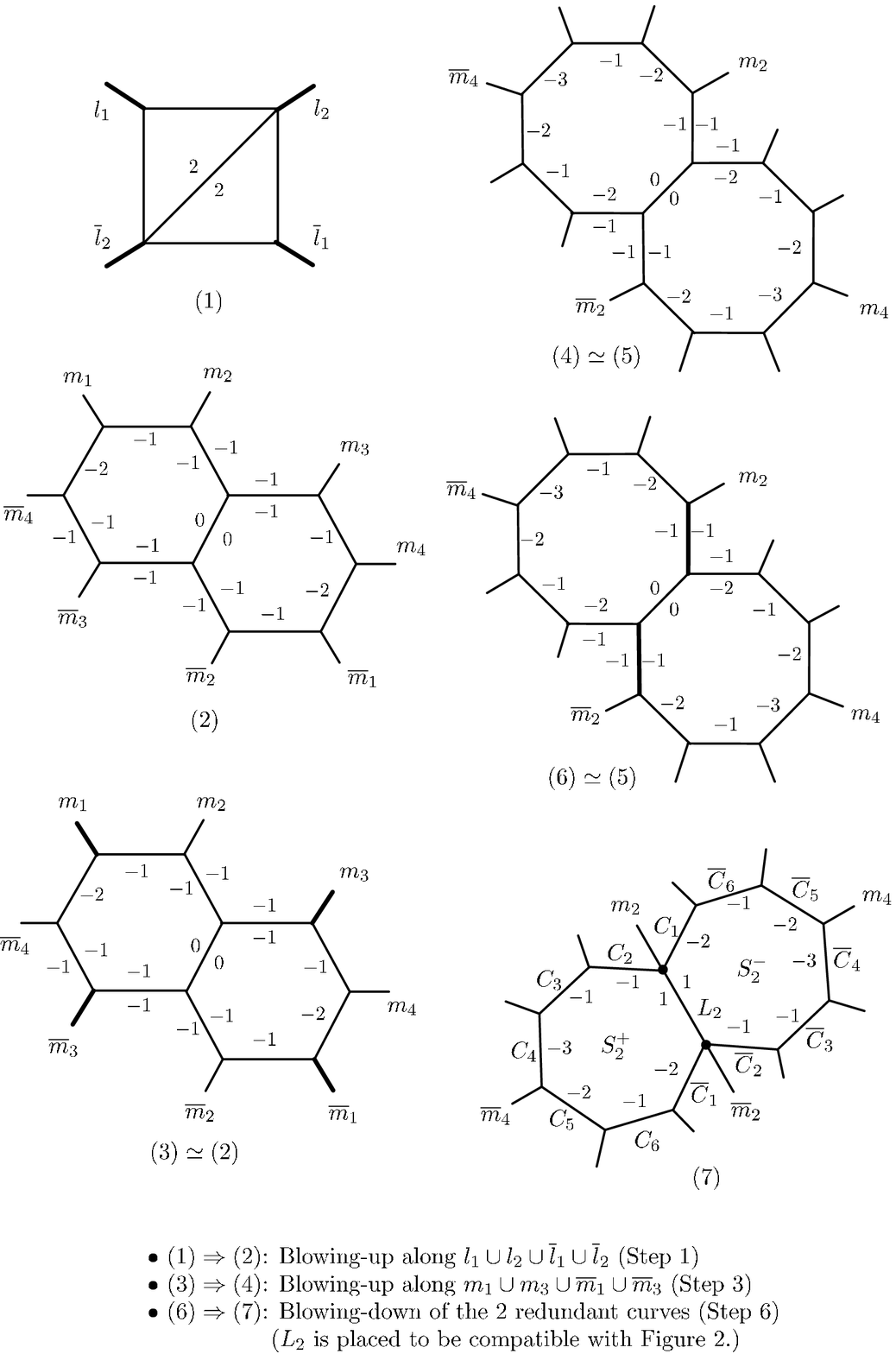}
\caption{The transformations for the reducible fiber over $\lambda_2$}
\label{fig-bim2}
\end{figure}

\begin{figure}
\includegraphics{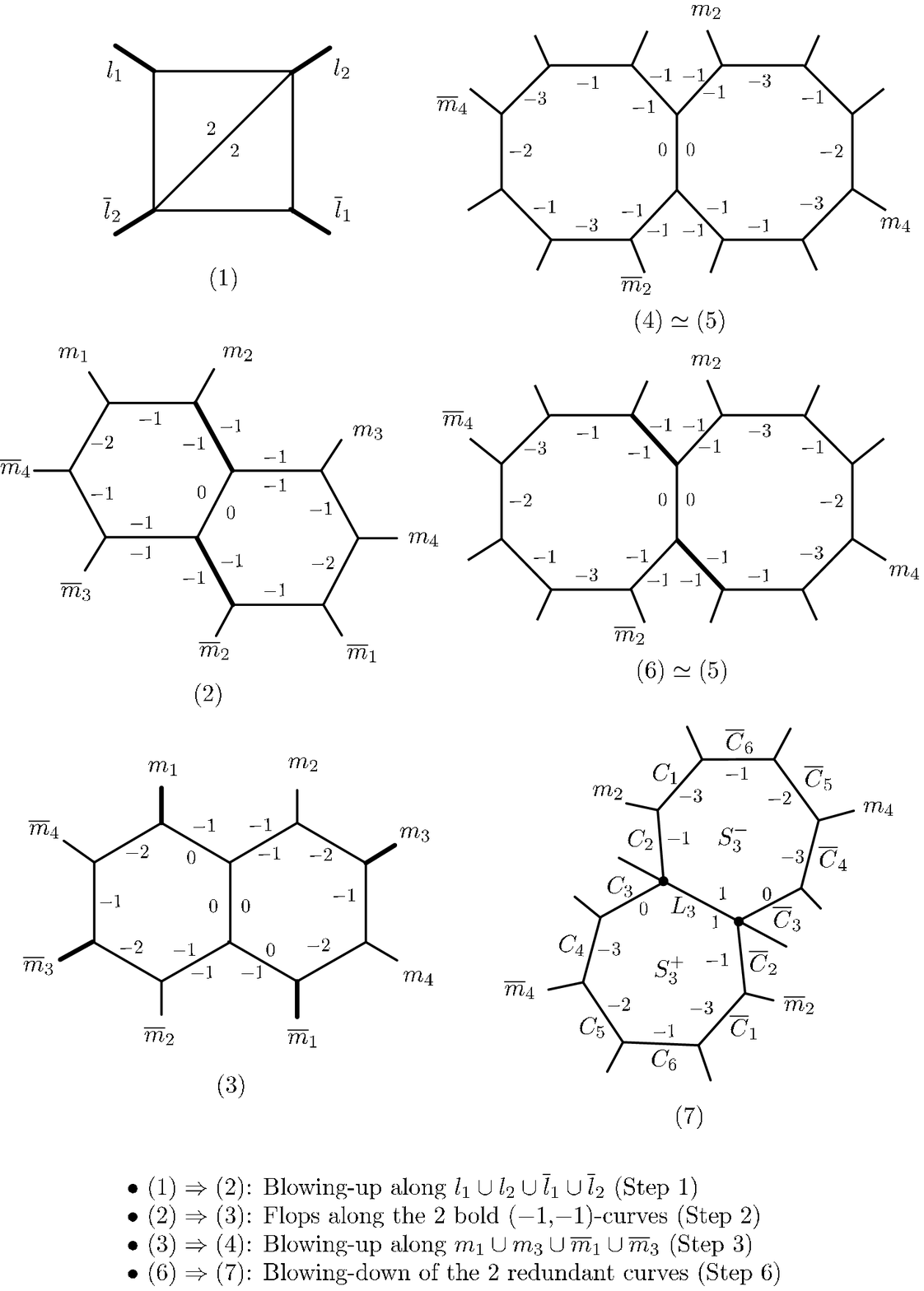}
\caption{The transformations for the reducible fiber over $\lambda_3$}
\label{fig-bim3}
\end{figure}

\begin{figure}
\includegraphics{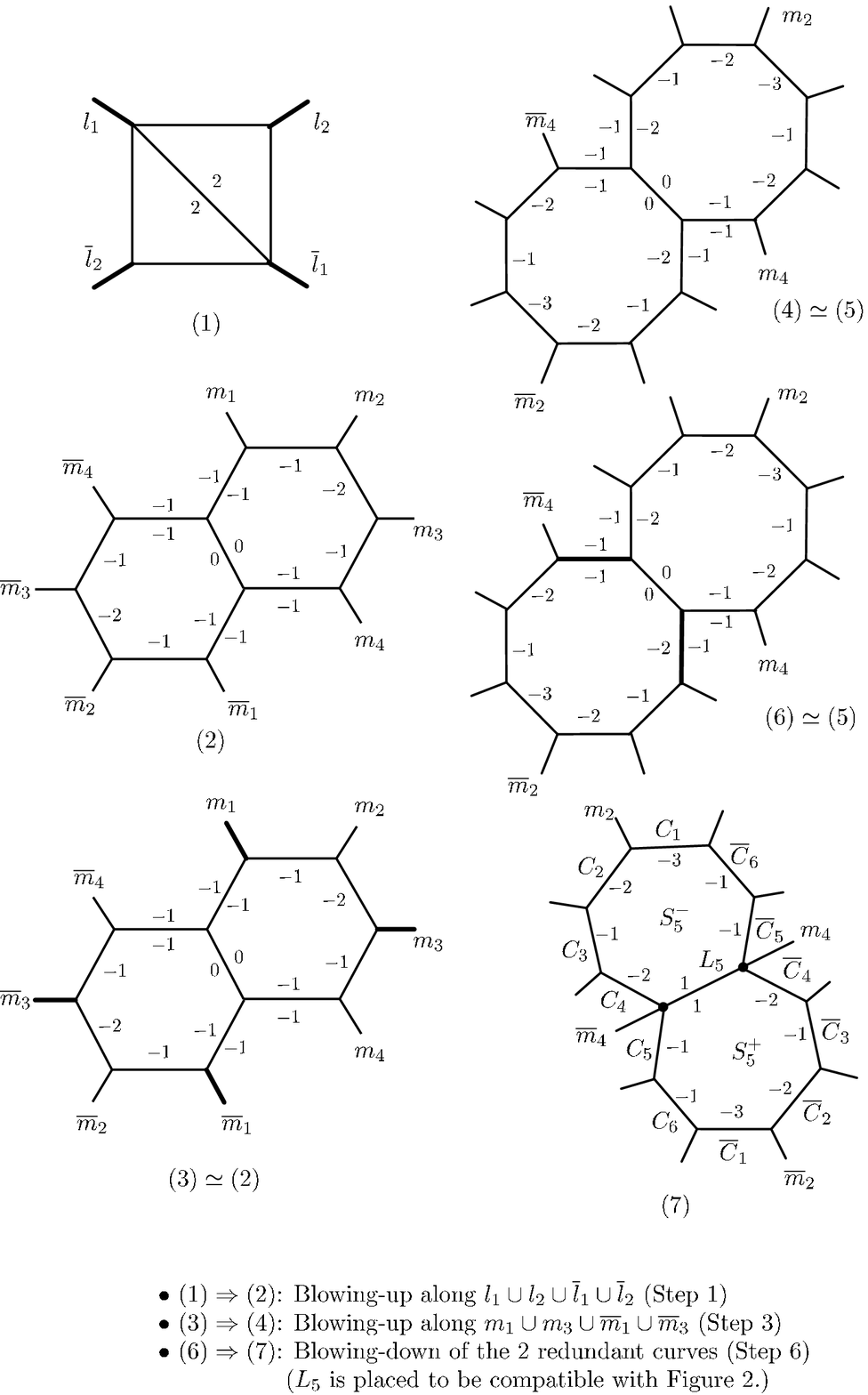}
\caption{The transformations for the reducible fiber over $\lambda_5$}
\label{fig-bim5}
\end{figure}

\begin{figure}
\includegraphics{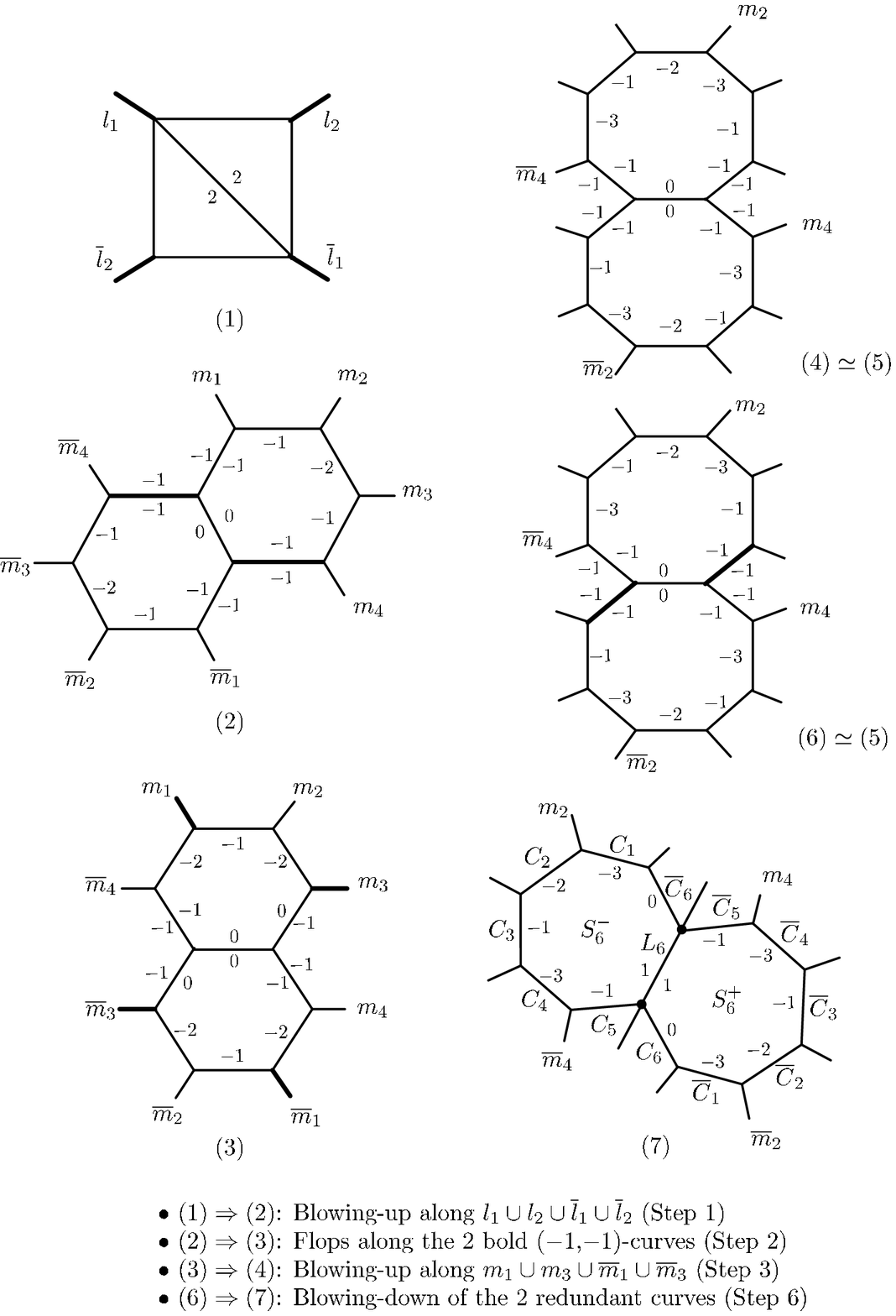}
\caption{The transformations for the reducible fiber over $\lambda_6$}
\label{fig-bim6}
\end{figure}

\noindent
(b)
If $i=2,3,5,6$, $p_1^{-1}(\lambda_i)$ is a reducible toric surface consisting of 2 quadratic cones ((1) of Figures \ref{fig-bim2},\ref{fig-bim3},\ref{fig-bim5},\ref{fig-bim6}).
These 2 cones form a conjugate pair, and their  intersection is a smooth real conic.
(These 4 conics will become $G$-invariant twistor lines.)

\begin{figure}
\includegraphics{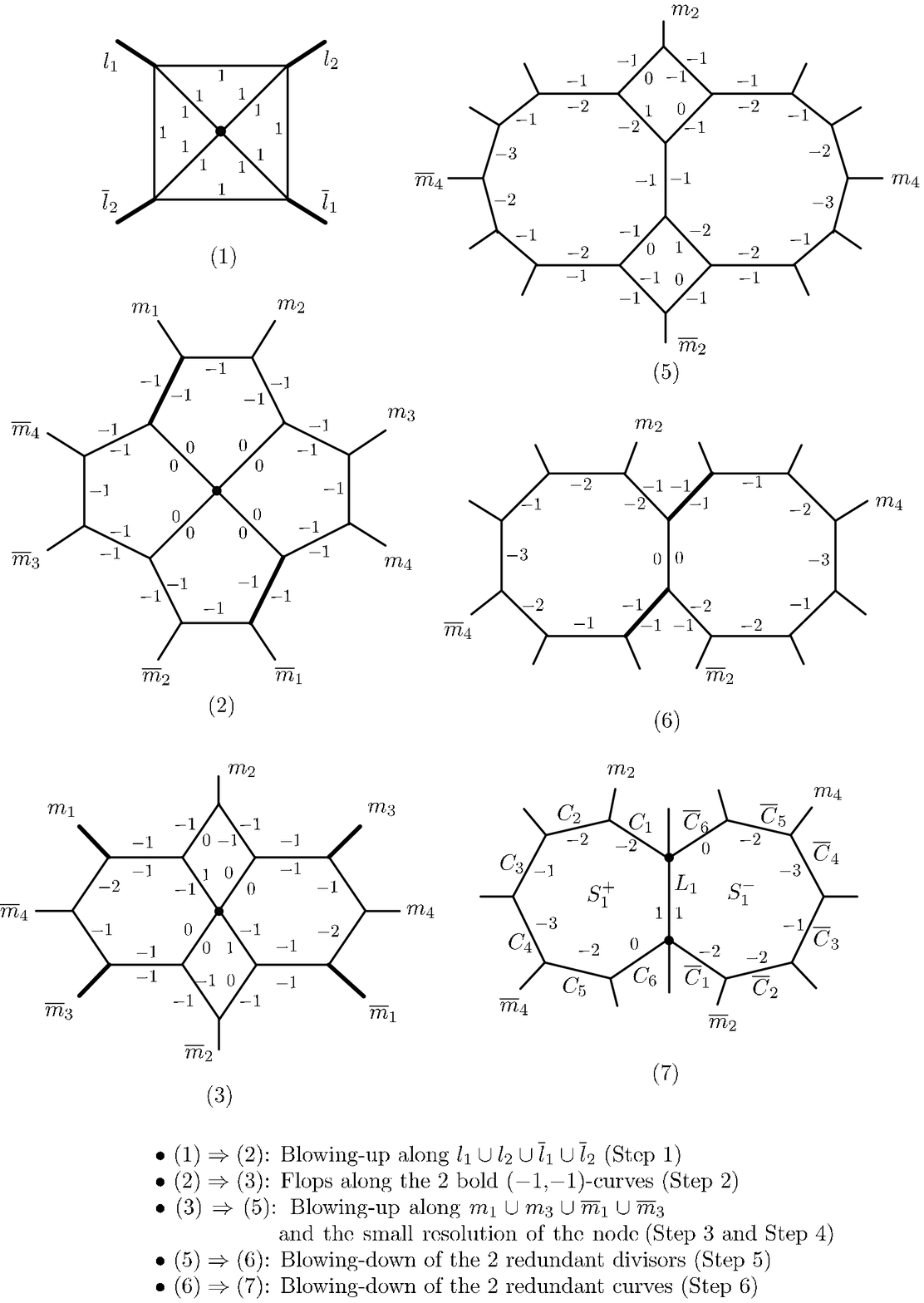}
\caption{The transformations for the reducible fiber over $\lambda_1$}
\label{fig-bim1}
\end{figure}

\begin{figure}
\includegraphics{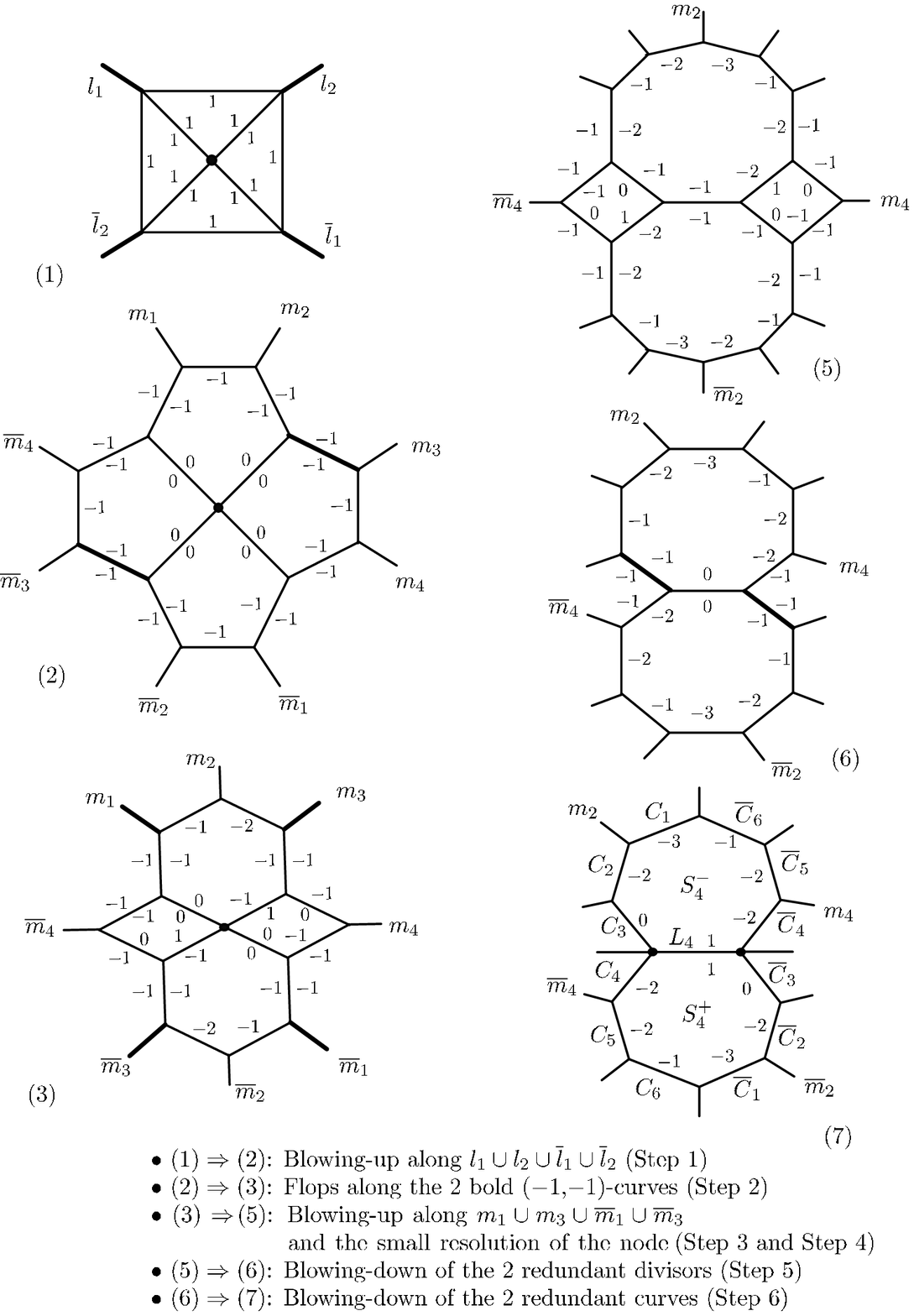}
\caption{The transformations for the reducible fiber over $\lambda_4$}
\label{fig-bim4}
\end{figure}

\noindent
(c)
If $i=1,4$, $p_1^{-1}(\lambda_i)$ is a reducible toric surface consisting of 4 projective planes ((1) of Figures \ref{fig-bim1},\ref{fig-bim4}).
These 4 planes share a unique point and it is an ordinary double point of $X_1$.
(These 2 points will yield the remaining two $G$-invariant twistor lines.)
The singular locus of $X_1$ consists of  $l_1\cup l_2\cup\ol{l}_1\cup \ol{l}_2$ and these 2 ordinary double points.
Thus the singular locus of $X_1$ coincides with the set of $G$-fixed points on $X_1$.

\vsp
In the following   we apply a succession of blowing-ups and blowing-downs to $X_1$.
The goal is to obtain a fiber space $\hat Z\ra\Lambda=\mathbf P^1$ of toric surfaces, which will be shown to be isomorphic to a blowup of the twistor space of a Joyce metric.
In every steps blowing-ups and downs are performed in such a way that they preserve the $G$-action and the real structure. 
So we do not mention about the $G$-action and the real structure in each step.

\vsp
\noindent
$\bullet$ {\bf{Step 1}} (Resolution of the 1-dimensional singular locus of $X_1$).
Let $X_2\ra X_1$ be the blowing-up along singularities $l_1\cup l_2\cup\ol{l}_1\cup \ol{l}_2$
and $p_2:X_2\ra\Lambda$ the natural projection.
On the exceptional divisors of $X_2\ra X_1$ there are 8 distinguished sections $m_1,\cdots,m_4$ and $\ol{m}_1,\cdots,\ol{m}_4$ of $p_2$  whose points are $G$-fixed ((2) of Figures \ref{fig-bim0}--\ref{fig-bim6}).
In the subsequent construction, these sections become important in that they enable us to keep track of mutual relationship between different reducible fibers.
Under the blowing-up $X_2\ra X_1$, each fiber gets the following effect:

\noindent
(a) On irreducible fibers of $p_1:X_1\ra\Lambda$, 
$X_2\ra X_1$ gives the minimal resolution of all the singularities therein
($(1)\Ra (2)$ of Figure \ref{fig-bim0}).

\noindent
(b) If $i=2,3,5,6$, the intersection point $p_1^{-1}(\lambda_i)\cap l_j$ and $p_1^{-1}(\lambda_i)\cap \ol{l}_j$
 ($j=1,2$) is either the vertex of the 2 cones or a smooth $G$-fixed point lying on the intersection conic of the cones.
For the former intersection points, $X_2\ra X_1$ gives the minimal resolution, whereas for the latter intersection points, it gives usual blow-ups (at smooth points) of both cones.
Thus for each $i=2,3,5,6$ the restriction $p_2^{-1}(\lambda_i)\ra p_1^{-1}(\lambda_i)$ of $X_2\ra X_1$ has 6 exceptional curves, two of which become $(-2)$-curves (which are the exceptional curves over the 2 vertices) in $p_2^{-1}(\lambda_i)$, and four of which are $(-1)$-curves in the strict transforms of the cones.
Moreover, the latter 4 curves  are $(-1,-1)$-curves in $X_2$
((2) of Figures \ref{fig-bim2},\ref{fig-bim3},\ref{fig-bim5},\ref{fig-bim6}).

\noindent
(c) If $i=1,4$, all the intersection points $p_1^{-1}(\lambda_i)\cap l_j$ and $p_1^{-1}(\lambda_i)\cap \ol{l}_j$
 ($j=1,2$)  are smooth points of two planes. 
So the exceptional locus of $p_2^{-1}(\lambda_i)\ra p_1^{-1}(\lambda_i)$ consists of 8 rational curves. (Each planes are blown-up twice.) 
All these exceptional curves are $(-1,-1)$-curves in $X_2$, and they intersect precisely one of the sections  $m_j$ and $\ol{m}_j$, $1\le j\le 4$.
((2) of Figures \ref{fig-bim1},\ref{fig-bim4})

\vsp
\noindent
$\bullet$ {\bf{Step 2}} (Flops inside some reducible fibers).
For each $i=1,3,4,6$ (not $2,3,5,6$) we will choose a conjugate pair of  $(-1,-1)$-curves in $p_2^{-1}(\lambda_i)$ among those obtained in Step 1.
We give these pairs by specifying which sections among $m_1,\cdots,m_4$ and $\ol{m}_1,\cdots,\ol{m}_4$ they intersect respectively:
 for $i=1$, we choose  a pair of $(-1,-1)$-curves intersecting $m_1$ and $\ol{m}_1$ respectively;
for $i=4$, we choose a pair of $(-1,-1)$-curves intersecting $m_3$ and $\ol{m}_3$ respectively;
for $i=3$, we choose a pair of $(-1,-1)$-curves intersecting $m_2$ and $\ol{m}_2$ respectively;
for $i=6$, we choose a pair of $(-1,-1)$-curves intersecting $m_4$ and $\ol{m}_4$ respectively.
Then we take flops along these 8 curves
($(2)\Ra(3)$ of Figures \ref{fig-bim1},\ref{fig-bim3},\ref{fig-bim4},\ref{fig-bim6}).
Let $X_3$ be the resulting space and $p_3:X_3\ra\Lambda$ the natural projection.
For the strict transforms of the sections $m_j$ and $\ol{m}_j$, $1\le j\le 4$ we use the same notations $m_j$ and $\ol{m}_j$ respectively.

\vsp
\noindent
$\bullet$ {\bf{Step 3}} (Further blowing-ups along 4 sections).
Let $X_4\ra X_3$ be the blowing-up along $m_1\cup m_3\cup \ol{m}_1\cup \ol{m}_3$ and $p_4:X_4\ra\Lambda$ the projection
($(3)\Ra(4)$ of the Figures \ref{fig-bim0}--\ref{fig-bim6}; in Figures \ref{fig-bim0} and \ref{fig-bim4}, the pictures of  `(4)' are not included because of the lack of space.)
Then irreducible fibers of $p_3$ are blown-up at 4 smooth points and consequently they become a (smooth) toric surface with $K^2=0$. (We do not touch these smooth fibers any further.)
For $i=1,4$,  among 4 irreducible components of $p_3^{-1}(\lambda_i)$, only one conjugate pair of the components get effects by $X_4\ra X_3$ and each of these 2 components are blown-up twice.
Another pair of the components of $p_3^{-1}(\lambda_i)$ are untouched through $X_4\ra X_3$.
So we temporary call these 2 components `the untouched components'.
(In (3) and (5) of Figures \ref{fig-bim1} and \ref{fig-bim4} these components are represented by the 2 squares.)
It is readily seen that these untouched components are isomorphic to $\Sigma_1$, a one point blown-up of $\mathbf P^2$.
On the other hand, for $i=2,3,5,6$, each irreducible components of $p_3^{-1}(\lambda_i)$ are blown-up at 2 points through $X_4\ra X_3$.

\vsp
\noindent
$\bullet$ {\bf{Step 4}} (Small resolutions of the ordinary double points).
The two ordinary double points of $X_1$ get no effect through Step 1--3, and $X_4$ has the corresponding ordinary double points which are on $p_4^{-1}(\lambda_1)$ and 
$p_4^{-1}(\lambda_4)$ respectively.
Let $X_5\ra X_4$ be their small resolutions which do {\em{not}} blow-up the `untouched components' in Step 3, and $p_5:X_5\ra\Lambda$ the projection.
This condition uniquely determine the small resolutions.
(The 2 exceptional curves of $X_5\ra X_4$ will become $G$-invariant twistor lines.)
Consequently $X_5$ becomes non-singular and the normal bundles of the strict transforms of the untouched components into $X_5$ has degree $-1$ along fibers of $\Sigma_1\ra\mathbf P^1$ ((5) of Figures  \ref{fig-bim1} and \ref{fig-bim4}).

\vsp
\noindent
$\bullet$ {\bf{Step 5}} (Blowing-downs of redundant divisors).
From the last remark in Step 4, the `untouched components' over $\lambda_1$ and $\lambda_4$ can be blown-down along a projection $\Sigma_1\ra\mathbf P^1$.
Let $X_5\ra X_6$ be the blowing-down of every 4 untouched components and $p_6:X_6\ra\Lambda$ the projection ($(5)\Ra (6)$ of Figures  \ref{fig-bim1} and \ref{fig-bim4}).
$X_6$ is still non-singular.
Another conjugate pair of the components over $\lambda_1$ and $\lambda_4$ get no effect through $X_5\ra X_6$.
Now all reducible fibers  $p_6^{-1}(\lambda_i),1\le i\le 6$, consist of 2 irreducible components which are conjugate of each other ((6) of Figures \ref{fig-bim1}--\ref{fig-bim6}).
The intersection of these components are real smooth rational curves in $p^{-1}_6(\lambda_i)$, which we write $L_i$, $1\le i\le 6$.
($L_i$ will become $G$-invariant twistor lines.)

\vsp
\noindent
$\bullet$ {\bf{Step 6}} (Contraction of redundant curves into ordinary double points).
For each $1\le i\le 6$, there is a conjugate pair of $(-1,-1)$-curves inside $p_6^{-1}(\lambda_i)$ which intersect $L_i$ (the bold curves in (6) of Figures \ref{fig-bim1}--\ref{fig-bim6}).
Let $X_6\ra \hat Z$ be the contraction of these twelve $(-1,-1)$-curves into ordinary double points and $f:\hat Z\ra\Lambda$ the natural projection.
Thus $\hat Z$ has twelve ordinary double points
(dotted points in (7) of Figures \ref{fig-bim1}--\ref{fig-bim6}).

\begin{rmk}
{\rm
As one may notice, the operations displayed in Figures \ref{fig-bim1}--\ref{fig-bim6} are not independent.
Each pictures in Figure \ref{fig-bim4} are obtained from those in Figures \ref{fig-bim1} by rotating $90^{\circ}$ clockwisely, while the place of all sections $l_j$ and $m_j$ are fixed.
There is the same relationship between Figure \ref{fig-bim5} and Figure \ref{fig-bim2}, and Figure \ref{fig-bim6} and Figure \ref{fig-bim3}
}\end{rmk}

\vsp
Now we have eventually reached the main result which provides an explicit construction of the twistor space of the Joyce metrics in problem: 

\begin{thm}
Let $g$ be a Joyce's self-dual metric on $4\mathbf P^2$ whose $K$-action is type I.
Let $\{\lambda_1<\lambda_2<\cdots<\lambda_6\}$ be the conformal invariant of\, $g$.
Let $Z$ be the twistor space of\, $g$, $X$ its anticanonical model which is the complete intersection in $\mathbf{P}^6$ defined by \eqref{732}, and $\hat Z$ the 3-fold  explicitly constructed from $X$ by Step 1-- Step 6 above. 
Let $\tilde Z\ra Z$ be the blowing-up  of $Z$ along the cycle $C$.
Then there is an isomorphism between $\tilde Z$ and $\hat Z$ which commutes with the real structues.
\end{thm}

\noindent Proof. 
We have a natural morphism $\tilde f:\tilde Z\ra\mathbf P^1$ whose fibers are members of the pencil $|(-1/2)K_Z|$ which are toric surfaces.
The structure  of all fibers of $\tilde f$ can be explicitly determined from the given $K$-action on $4\mathbf P^2$.
On the other hand, by our construction there is  a natural morphism $f:\hat Z\ra \Lambda=\mathbf P^1$ whose fibers are also toric surfaces.
Since all the operations in Steps 1-- 6 are explicit, we can verify that all fibers of $f$ and $\tilde f$ are $G$-equivariantly isomorphic, there is no real point on $\hat Z$, and that the real structure on $\hat Z$ exchanges $C_i$ and $\ol{C}_i$ contained in real fibers.
Then by \cite[Section 8]{F00}, 
in order to obtain an isomorphism between $\hat Z$ and $\tilde Z$ which commutes with the real structures, it suffices to prove the existence of a
particular curve in $\hat Z$ which plays a vital role  in  \cite[Section 8]{F00} to produce an isomorphism between two distinct fiber spaces of toric surfaces whose fibers are mutually isomorphic.
Namely, it suffices to show that there exists a smooth rational curve $\hat L$ in $\hat Z$ fulfilling the following properties:
(i) $\hat L$ is a smooth real rational curve,
(ii) the restriction of the projection $f:\hat Z\ra\Lambda$ onto $\hat L$ is 2 to 1 whose branch points become a pair of conjugate points in $\Lambda\backslash\Lambda^{\sigma}$.
(iii) $\hat L$ does not intersect the strict transform of $\nu_1^{-1}(B)$ into $\hat Z$, where $\nu_1:X_1\ra X$ is the blowup along $B$ as before.

To show the existence of such a curve $\hat L$,  take  a general twistor line in $Z$ in the sense of \cite[(6.11)]{F00}.
Namely let $L$ be a twistor line  which does not intersect $C$.
Then $L$ is contained in $U$ (= the union of 2-dimensional $G$-orbits introduced in the final part of Section 2.)
Hence by Lemma \ref{lemma-32} (i), $\Phi(L)$  is a real smooth rational curves contained in $X\backslash B$.
Therefore $\nu_1^{-1}(\Phi(L))$ does not intersect $\nu_1^{-1}(B)$.
Then we define $\hat L\subset \hat Z$ to be the strict transform of $\nu_1^{-1}(B)$ under the birational transformations in Step 1--6.
Then since $L\cdot(-1/2)K_Z=2$, the restriction of the projection $\hat Z\ra\Lambda$ onto $\hat L$ is 2 to 1 and its ramification locus can supposed to be in $\Lambda\backslash\Lambda^{\sigma}$ by choosing sufficiently general $L$.
Thus $\hat L$ satisfies (ii).
To verify (i) and (iii) we have to see that Step 1-- Step 6 do not change an open neighborhood of $\nu_1^{-1}(\Phi(L))$. 
This is obvious except in Step 5.
In Step 5 we have blown-down  the `redundant' divisors (namely the `untouched components'), and we have to check that the strict transform of $\nu_1^{-1}(\Phi(L))$ into $X_4$ does not go through these  divisors.
To verify this, let $D^+_1+D^-_1$ and $D^+_4+D^-_4$ be the redundant divisors in   $p_4^{-1}(\lambda_1)$ and  $p_
4^{-1}(\lambda_4)$ respectively.
Then by our choice of the pair of $(-1,-1)$-curves in Step 2,  $D^+_1+D^-_1$ are put in  a skew positions in relation to $D^+_4+D^-_4$; namely recalling that $D_1^{\pm}$ and $D_4^{\pm}$ correspond (through the birational changes) to planes in $X\subset \mathbf P^6$, 
the 4 planes   intersect $\mathbf P^3_{\infty}$ along 4 different lines.
On the other hand, as showed in Lemma \ref{lemma-11}, the set $\{\Phi(S_i^+)=\mathbf P^{\vee}W_i^+, \Phi(S_i^-)=\mathbf P^{\vee}W_i^-\set i=1,4\}$ becomes 4 planes (in $\mathbf P^6$).
Further by Lemma \ref{lemma-32} (ii) 
the intersection of these 4 planes with $\mathbf P^3_{\infty}$ is 4 different lines.
These imply that $\Phi(L)$ does not intersect the redundant planes in $X$ at all, or else
intersect all of the redundant planes.
If the latter situation actually happens, we exchange the role of $\lambda_1$, $\lambda_2$, $\lambda_3$, $\lambda_4$, $\lambda_5$ and $\lambda_6$, by $\lambda_4,\lambda_3,\lambda_2,\lambda_1, \lambda_6$ and $\lambda_5$ respectively, in every steps.
Then the pairs redundant planes are exchanged and hence $\Phi(L)$ becomes to be not intersecting redundant planes.
Thus $\hat L$ is shown to satisfy (i)--(iii) and hence it follows that $\hat Z$ and $\tilde Z$ are isomorphic.
\proofend



\small
\vspace{13mm}
\hspace{5.5cm}
$\begin{array}{l}
\mbox{Department of Mathematics}\\
\mbox{Graduate School of Science and Engineering}\\
\mbox{Tokyo Institute of Technology}\\
\mbox{2-12-1, O-okayama, Meguro, 152-8551, JAPAN}\\
\mbox{{\tt {honda@math.titech.ac.jp}}}
\end{array}$

\newpage

\end{document}